\def\DATE{Nov.\ 24, 2012}
\def\DATE{\relax}
\magnification=1100
\baselineskip=12.72pt
\voffset=.75in
\hoffset=.9in
\hsize=4.1in
\newdimen\hsizeGlobal
\hsizeGlobal=4.1in%
\vsize=7.05in
\parindent=.166666in
\pretolerance=500 \tolerance=1000 \brokenpenalty=5000

\footline={\vbox{\hsize=\hsizeGlobal\hfill{\rm\the\pageno}\hfill\llap{\sevenrm\DATE}}\hss}

\def\note#1{%
  \hfuzz=50pt%
  \vadjust{%
    \setbox1=\vtop{%
      \hsize 3cm\parindent=0pt\eightpoints\baselineskip=9pt%
      \rightskip=4mm plus 4mm\raggedright#1\hss%
      }%
    \hbox{\kern-4cm\smash{\box1}\hss\par}%
    }%
  \hfuzz=0pt
  }
\def\note#1{\relax}

\def\anote#1#2#3{\smash{\kern#1in{\raise#2in\hbox{#3}}}%
  \nointerlineskip}     

\newcount\equanumber
\equanumber=0
\newcount\sectionnumber
\sectionnumber=0
\newcount\subsectionnumber
\subsectionnumber=0
\newcount\snumber  
\snumber=0

\def\section#1{%
  \subsectionnumber=0%
  \snumber=0%
  \equanumber=0%
  \advance\sectionnumber by 1%
  \noindent{\bf \the\sectionnumber .~#1.~}%
}%
\def\subsection#1{%
  \advance\subsectionnumber by 1%
  \snumber=0%
  \equanumber=0%
  \noindent{\bf \the\sectionnumber .\the\subsectionnumber .~#1.~}%
}%
\def\prevs{\the\sectionnumber .\the\subsectionnumber .\the\snumber }
\long\def\Definition#1{%
  \global\advance\snumber by 1%
  \bigskip
  \noindent{\bf Definition~\prevs .}%
  \quad{\it#1}%
}

\long\def\Corollary#1{%
  \global\advance\snumber by 1%
  \bigskip
  \noindent{\bf Corollary~\prevs .}%
  \quad{\it#1}%
}%
\long\def\Lemma#1{%
  \global\advance\snumber by 1%
  \bigskip
  \noindent{\bf Lemma~\prevs .}%
  \quad{\it#1}%
}%
\def\Proof{\noindent{\bf Proof.~}}
\long\def\Proposition#1{%
  \advance\snumber by 1%
  \bigskip
  \noindent{\bf Proposition~\prevs .}%
  \quad{\it#1}%
}%
\long\def\Remark#1{%
  \bigskip
  \noindent{\bf Remark.~}#1%
}%
\long\def\remark#1{%
  \advance\snumber by1%
  \bigskip
  \noindent{\bf Remark~\prevs .}\quad#1%
}%
\long\def\Theorem#1{%
  \advance\snumber by 1%
  \bigskip
  \noindent{\bf Theorem~\prevs .}%
  \quad{\it#1}%
}%
\long\def\Statement#1{%
  \advance\snumber by 1%
  \bigskip
  \noindent{\bf Statement~\prevs .}%
  \quad{\it#1}%
}%
\def\ifundefined#1{\expandafter\ifx\csname#1\endcsname\relax}
\def\labeldef#1{\global\expandafter\edef\csname#1\endcsname{\prevs}}
\def\labelref#1{\expandafter\csname#1\endcsname}
\def\label#1{\ifundefined{#1}\labeldef{#1}\note{$<$#1$>$}\else\labelref{#1}\fi}

\def\preveq{(\the\sectionnumber .\the\subsectionnumber .\the\equanumber)}
\def\neq{\global\advance\equanumber by 1\eqno{\preveq}}

\def\ifundefined#1{\expandafter\ifx\csname#1\endcsname\relax}

\def\equadef#1{\global\advance\equanumber by 1%
  \global\expandafter\edef\csname#1\endcsname{\preveq}%
  \setbox1=\hbox{\rm\hskip .1in[#1]}\dp1=0pt\ht1=0pt\wd1=0pt%
  \preveq\box1}
\def\equadef#1{\global\advance\equanumber by 1%
  \global\expandafter\edef\csname#1\endcsname{\preveq}%
  \preveq}

\def\equaref#1{\expandafter\csname#1\endcsname}

\def\equa#1{%
  \ifundefined{#1}%
    \equadef{#1}%
  \else\equaref{#1}\fi}

\font\eightrm=cmr8%
\font\sixrm=cmr6%

\font\eightsl=cmsl8%

\font\eightbf=cmb8%

\font\eighti=cmmi8%
\font\sixi=cmmi6%

\font\eightsy=cmsy8%
\font\sixsy=cmsy6%

\font\eightex=cmex8%
\font\sixex=cmex6%
\font\fiveex=cmex5%

\font\eightit=cmti8%

\font\eighttt=cmtt8%

\font\tenbb=msbm10%
\font\eightbb=msbm8%
\font\sevenbb=msbm7%
\font\sixbb=msbm6%
\font\fivebb=msbm5%
\newfam\bbfam  \textfont\bbfam=\tenbb  \scriptfont\bbfam=\sevenbb  \scriptscriptfont\bbfam=\fivebb%

\font\tenbbm=bbm10

\font\tencmssi=cmssi10%
\font\sevencmssi=cmssi7%
\font\fivecmssi=cmssi5%
\newfam\ssfam  \textfont\ssfam=\tencmssi  \scriptfont\ssfam=\sevencmssi  \scriptscriptfont\ssfam=\fivecmssi%

\font\tenfrak=cmfrak10%
\font\eightfrak=cmfrak8%
\font\sevenfrak=cmfrak7%
\font\sixfrak=cmfrak6%
\font\fivefrak=cmfrak5%
\newfam\frakfam  \textfont\frakfam=\tenfrak  \scriptfont\frakfam=\sevenfrak  \scriptscriptfont\frakfam=\fivefrak%
\def\frak{\fam\frakfam\tenfrak}%

\font\tenmsam=msam10%
\font\eightmsam=msam8%
\font\sevenmsam=msam7%
\font\sixmsam=msam6%
\font\fivemsam=msam5%

\def\bb{\fam\bbfam\tenbb}%

\def\hexdigit#1{\ifnum#1<10 \number#1\else%
  \ifnum#1=10 A\else\ifnum#1=11 B\else\ifnum#1=12 C\else%
  \ifnum#1=13 D\else\ifnum#1=14 E\else\ifnum#1=15 F\fi%
  \fi\fi\fi\fi\fi\fi}
\newfam\msamfam  \textfont\msamfam=\tenmsam  \scriptfont\msamfam=\sevenmsam  \scriptscriptfont\msamfam=\fivemsam%
\def\msam{\msamfam\tenmsam}%
\mathchardef\leq"3\hexdigit\msamfam 36%
\mathchardef\geq"3\hexdigit\msamfam 3E%

\font\tentt=cmtt11%
\font\seventt=cmtt9%
\textfont\ttfam=\tentt
\scriptfont7=\seventt%
\def\tt{\fam\ttfam\tentt}%

\def\eightpoints{%
\def\rm{\fam0\eightrm}%
\textfont0=\eightrm   \scriptfont0=\sixrm   \scriptscriptfont0=\fiverm%
\textfont1=\eighti    \scriptfont1=\sixi    \scriptscriptfont1=\fivei%
\textfont2=\eightsy   \scriptfont2=\sixsy   \scriptscriptfont2=\fivesy%
\textfont3=\eightex   \scriptfont3=\sixex   \scriptscriptfont3=\fiveex%
\textfont\itfam=\eightit  \def\it{\fam\itfam\eightit}%
\textfont\slfam=\eightsl  \def\sl{\fam\slfam\eightsl}%
\textfont\ttfam=\eighttt  \def\tt{\fam\ttfam\eighttt}%
\textfont\bffam=\eightbf  \def\bf{\fam\bffam\eightbf}%

\textfont\frakfam=\eightfrak  \scriptfont\frakfam=\sixfrak \scriptscriptfont\frakfam=\fivefrak  \def\frak{\fam\frakfam\eightfrak}%
\textfont\bbfam=\eightbb      \scriptfont\bbfam=\sixbb     \scriptscriptfont\bbfam=\fivebb      \def\bb{\fam\bbfam\eightbb}%
\textfont\msamfam=\eightmsam  \scriptfont\msamfam=\sixmsam \scriptscriptfont\msamfam=\fivemsam  \def\msam{\msamfam\eightmsam}

\rm%
}

\def\poorBold#1{\setbox1=\hbox{#1}\wd1=0pt\copy1\hskip.25pt\box1\hskip .25pt#1}

\mathchardef\lsim"3\hexdigit\msamfam 2E%
\mathchardef\gsim"3\hexdigit\msamfam 26%

\def\Ai{{\rm Ai}}

\def\d{\,{\rm d}}

\def\ds{\displaystyle}
\long\def\DoNotPrint#1{\relax}
\def\e{{\rm e}}

\def\fixedref#1{#1\note{fixedref$\{$#1$\}$}}

\def\Id{{\rm Id}}

\def\lstar#1{{}^*\hskip -2.8pt#1}

\def\qed{{\vrule height .9ex width .8ex depth -.1ex}}
\def\sS{\scriptstyle}

\def\boc{\note{{\bf BoC}\hskip-11pt\setbox1=\hbox{$\Bigg\downarrow$}%
         \dp1=0pt\ht1=0pt\ht1=0pt\leavevmode\raise -20pt\box1}}
\def\eoc{\note{{\bf EoC}\hskip-11pt\setbox1=\hbox{$\Bigg\uparrow$}%
         \dp1=0pt\ht1=0pt\ht1=0pt\leavevmode\raise 20pt\box1}}

\def\One{\hbox{\tenbbm 1}}

\def\calC{{\cal C}}

\def\calT{{\cal T}}

\def\CC{{\bb C}}

\def\MM{{\bb M}\kern .4pt}
\def\NN{{\bb N}\kern .5pt}
\def\RR{{\bb R}}

\def\ZZ{{\bb Z}}

\pageno=1

\centerline{\bf \poorBold{$q$}-CATALAN BASES}
\centerline{\bf AND THEIR DUAL COEFFICIENTS}

\bigskip
 
\centerline{Ph.\ Barbe$^{(1)}$ and W.P.\ McCormick$^{(2)}$}
\centerline{${}^{(1)}$CNRS {\sevenrm(UMR {\eightrm 8088})}, ${}^{(2)}$University
 of Georgia}

 
{\narrower
\baselineskip=9pt\parindent=0pt\eightpoints

\bigskip

{\bf Abstract.} We define $q$-Catalan bases which are a generalization of
the $q$-polyomials $z^n(z,q)_n$. The determination of their dual bases
involves some $q$-power series termed dual coefficients. We show how these
dual coefficients occur in the solution of some equations with $q$-commuting
coefficients and solve an abstract $q$-Segner recursion. We study the 
connection between this theory and Garsia's (1981). The overall flavor of 
this work is to show how some properties of $q$-Catalan numbers are in fact
instances of much more general results on dual coefficients.

\bigskip


\bigskip
 

}

\bigskip

\def\prevs{\the\sectionnumber.\the\snumber }
\def\preveq{(\the\sectionnumber.\the\equanumber)}

\section{Introduction}
The purpose of this paper is to revisit the old problem of determining
the dual bases of some bases of power series, show some new connections
with other problems and give a different perspective on some results
pertaining to analytic combinatorics and the theory of $q$-series. Numerous 
papers deal with describing dual bases, but to broadly fix the ideas, while 
this work is not properly on combinatorics, it is inspired by results in that
area, especially the works of Garsia (1981), Krattenthaler (1984) and 
F\"urlinger and Hofbauer (1985) related to $q$-Catalan numbers.

In this paper we identify two new objects, a class of formal power
series, which we call $q$-Catalan bases, and some $q$-polynomials or
power series, which we call dual coefficients. While these dual
coefficients are instrumental in determining the dual bases of
$q$-Catalan bases, they are interesting quantities to study for their
own sake, since we will see that they connect different problems, such
as finding power series expansions to the solutions to some equations
with $q$-commuting coefficients and solving some recursions or
functional equations of $q$-type. Moreover, in some special cases,
these dual coefficients have a combinatorial interpretation such as
Carlitz's (1972) $q$-Catalan numbers, and our result leads to a new
characterization of these numbers in terms of power series solution of
a quadratic equation with $q$-commuting coefficients.

The organization of this paper is as follows. Section 2 contains the main 
definitions, in particular that of $q$-Catalan basis and dual coefficients,
and the main result of this section is a description of the dual bases of
$q$-Catalan bases. In section 3 we examine further properties of the dual
coefficients. Because the results of sections 2 and 3 are so intimately
related, section 4 gathers their proofs. In section 5 we consider some
particular $q$-Catalan bases; we show how some previous results by
Garsia (1981) and Krattenthaler (1988) can be viewed in this setting, and
how the Rogers-Ramanujan continuous fraction (see Flajolet and Sedgewick, 
2009; example V.9) can be extended
in this setting. Finally section 6 discusses further specialization related
to the combinatorics of trees and lattice paths.

\bigskip

\noindent{\bf Notation.} Following the custom in $q$-series (see e.g.\ Andrews,
Askey, Roy, 2000), we write for any $n$ positive 
$$
  (z,q)_n=\prod_{0\leq j<n} (1-zq^j) \, ,
$$ 
with $(z,q)_0=1$ and $(z,q)_\infty=\prod_{j\geq 0} (1-zq^j)$. 

When indexing quantities by pairs $(i,j)$, such as $e_{(i,j)}$ we tend to 
drop the pair notation and write instead $e_{i,j}$. We write $\CC[[z]]$ for the
set of all formal power series with complex coefficients.

Throughout the paper, we will consider formal power series. Conditions 
under which those formal power series are convergent ones are usually easy
to determine.

\bigskip

%
\def\prevs{\the\sectionnumber .\the\snumber }
\def\preveq{(\the\sectionnumber .\the\equanumber)}

\section{\poorBold{$q$}-Catalan bases and their dual forms}
Our main object of study is a special type of power series in two variables
and some related bases of the space of power series.

\Definition{\label{CatalanType}
  A power series $P(z,t)$ is a Catalan power series if there exists a power 
  series $\tilde P(z,t)$ such that 
  $$
    P(z,t)=t- z\tilde P(z,t)t^2 \, .
  $$
  We say that $\tilde P$ is associated to $P$.
}

\bigskip

While some nontrivial examples are developed in the fifth and sixth section, 
we will run a couple of trivial ones through this section in order to make 
things more concrete and explain the terminology.

\bigskip

\noindent{\bf Examples.} a) $P(z,t)=t$ is a Catalan power series with
$\tilde P(z,t)=0$.

\noindent b) $P(z,t)=t-zt^2$ is a Catalan power series with $\tilde P(z,t)=1$.

\bigskip

Viewing $P(z,t)$ as a power series in $t$ with coefficients in $\CC[[z]]$, 
it is a Catalan power series if the coefficient of $t$ is the power series
in $z$ which is constant and equal to $1$. If $P(z,t)$ and $Q(z,t)$ are
two Catalan power series, so is $\alpha P(z,t)+(1-\alpha)Q(z,t)$ for any 
real number $\alpha$, and so are $P(z,t)Q(z,t)/t$ and $t^2/P(z,t)$.

\bigskip

Next, we define what will turn out to be bases for the space of power series.

\Definition{\label{qCatalanType}
  A family of power series $\bigl(e_k(z,q)\bigr)_{k\geq 0}$ is a $q$-Catalan
  basis if there exists a Catalan power series $P$ such that for any 
  nonnegative integer $k$,
  $$
    {e_k(qz)\over e_k(z)}={P(z,q^k)\over P(z,1)} \, .
    \eqno{\equa{qCatalanTypeA}}
  $$
  We then say that $P$ is associated to $e_k$, or also that $e_k$ is
  associated to $P$.
}

\bigskip

\noindent{\bf Examples.} (continued) a) $e_k(z)=z^k$.

\noindent b) $e_k(z)=z^k(z,q)_k$.

\bigskip

Note that in both examples $e_k(z)$ is a polynomial in $z$, of order $k$.
Our first lemma shows that if $(e_k)$ is a Catalan basis then $e_k$ is of 
order $k$, and, consequently, that $(e_k)$ is a basis for the space of power 
series; hence the terminology.

\Lemma{\label{ekBasis}
  If $(e_k)$ is a $q$-Catalan basis, then each $e_k$ has order $k$.
}

\bigskip

To understand better the relationship between Catalan power series and
$q$-Catalan bases, we need the following definition.

\Definition{\label{canonical}
  A $q$-Catalan basis $(e_k)$ is normalized if $[z^k]e_k=1$.
}

\bigskip

A $q$-Catalan basis is normalized if its term of lowest degree has coefficient
$1$.

Among other things, our next result asserts that there is a bijection 
between Catalan power series and normalized $q$-Catalan bases.

\Lemma{\label{representation}
  (i) If $P$ is a Catalan power series, the unique normalized
  $q$-Catalan basis associated to $P$ is given by
  $$
    e_k(z)=z^k\prod_{j\in\NN} 
           {1-q^jz \tilde P(q^jz,1)\over 1-q^{k+j}z\tilde P(q^jz,q^k)} \, .
  $$

  \noindent (ii) A $q$-Catalan basis is associated with a unique Catalan
  power series.
}

\bigskip

Lemma \representation\ makes clear what the difference between $q$-Catalan
bases and Garsia's (1981) powers is: with the notation of 
Lemma \representation, Garsia's powers are obtained if $\tilde P(z,t)$ is
a power series in the product $zt$. In general, $q$-powers
of Hofbauer (1984), as extended by Krattenthaler (1984), are not $q$-Catalan
bases. Krattenthaler's (1988) results cover $q$-Catalan bases, but, not 
surprisingly, our less general assumptions entail for a more specific
theory which emphasizes less on the linear algebraic part of the theory
and more on its connection with other problems, and which is perhaps easier 
to apply in some cases. Our assumptions should be viewed as
in between those of Garsia (1981) and Krattenthaler (1988), though closer
to Garsia's (1981).

Having a basis, a natural question is how to decompose functions in this
basis. If $(e_n)$ is a basis, recall that its dual basis consists of the 
linear forms
$([e_n])$ such that $[e_n](e_k)=\delta_{k,n}$. Thus the question of decomposing
functions in a basis $(e_n)$ is really that of finding the dual basis.

Concerning the notation, note that when we consider the basis of
the monomials $(z^k)$ then $([z^k])$ are the corresponding dual forms, and
$[z^k]f(z)$ is the coefficient of $z^k$ in $f$. For power series of two 
variables $(z,t)$ say, then $[(z^it^j)]f(z,t)$ is the coefficient 
of $z^it^j$ in the power series $f(z,t)$.

In order for us to describe the dual bases of $q$-Catalan ones, we need
to define some functions of two variables. For this purpose, let
$$
  \e_1=(1,0)\qquad\hbox{and}\qquad \e_2=(0,1)
$$
be the canonical basis of $\RR^2$, and let $\langle\cdot,\cdot\rangle$ be
the usual inner product in $\RR^2$.\note{predual / prodromous / prefatory / elemental}

\Definition{\label{preDualBasis}
  Let $P$ be a Catalan power series. Its predual basis is the
  array $\bigl(\tilde e_i(z,t)\bigr)_{i\in\NN^2}$ of
  power series of two variables
  $$
    \tilde e_i(z,t)
    =z^{\langle i,\e_1\rangle} \prod_{0\leq j< \langle i,\e_2\rangle}
     P(q^jz,t)\, ,
  $$
  where it is agreed that a product over an empty index set is $1$.

  By extension, we also say that $(\tilde e_i)_{i\in\NN^2}$ is 
  a predual basis for any Catalan basis associated to $P$.
}

\bigskip

While the term predual is convenient and as we will see meaningful, it has no
connection with the algebraic notion of predual of a vector space.

\bigskip

\noindent{\bf Examples.} (continued) a) $\tilde e_i(z,t)
=z^{\langle i,\e_1\rangle}t^{\langle i,\e_2\rangle}$.

\noindent b) $\tilde e_i(z,t)=z^{\langle i,\e_1\rangle} t^{\langle i,\e_2\rangle}
(zt,q)_{\langle i,\e_2\rangle}=z^{\langle i,\e_1-\e_2\rangle} e_{\langle i,\e_2\rangle}(zt)$.

\bigskip

Given the convention that a product over an empty set is $1$, we have
$\tilde e_{(n,0)}=z^n$ for all nonnegative integers $n$.

Our next result asserts that a predual basis is indeed a basis.

\Lemma{\label{eTildeBasis}
   If $(e_k)_{k\in\NN}$ is a Catalan basis, its predual basis
   $\bigl(\tilde e_i(z,t)\bigr)_{i\in \NN^2}$ is a basis of the vector 
   space of power series in $(z,t)$.
}

\bigskip

As a consequence of Lemma \eTildeBasis, the following definition makes sense.

\Definition{\label{dualCoefficients}
  Given a Catalan power series $P$, or a Catalan basis $(e_k)_{k\in\NN}$ 
  associated to a Catalan
  power series $P$, its dual coefficients are the unique
  $(T_i)_{i\in\NN^2}$ defined by
  $$
    \sum_{i\in\NN^2} T_i \tilde e_i(z,t)=t \, .
    \eqno{\equa{dualCoefficientsEq}}
  $$
}

Note that since $(\tilde e_k)$ depends on $q$, the dual coefficients 
depend also on $q$, but not on $z$ or $t$.
A recursive way of calculating the dual coefficients is given in 
Theorem \fixedref{3.3} in the next section.

\bigskip

\noindent{\bf Examples.} (continued) a) $\sum_{i\in\NN^2} T_i 
z^{\langle i,\e_1\rangle} t^{\langle i,\e_2\rangle}=t$ forces $T_i=0$ if $i\not=(0,1)$
and $T_{0,1}=1$.

\noindent b) Multiplying both sides of \dualCoefficientsEq\ by $z$, we
obtain 
$$
  \sum_{i\in\NN^2} T_i z^{\langle i,\e_1-\e_2\rangle+1} 
  \e_{\langle i,\e_2\rangle}(zt)=zt \, .
  \eqno{\equa{exBa}}
$$ 
Since $(\e_n)$ is a basis and $e_0=1$, there exists
a unique sequence $(\calC_n)$ such that 
$$
  \sum_{n\geq 0} \calC_n e_{n+1}(z)=z \, .
  \eqno{\equa{exBb}}
$$ 
Then, for \exBa\ to hold, we need $T_i$ to vanish when 
$\langle i,\e_1-\e_2\rangle+1$
does not. Thus, $T_i=0$ if $i\not\in\{\, (n,n+1)\,:\, n\in\NN\,\}$ and for
$n$ in $\NN$, $T_{n,n+1}=\calC_n$. The sequence $(\calC_n)_{n\in\NN}$ satisfying
\exBb\ is known as Carlitz's (1972) $q$-Catalan numbers and does not have
a known nice closed form. It has been studied by Andrews (1975),
F\"urlinger and Hofbauer (1985) among others.

\bigskip

While the specific values of the dual coefficients depend on the Catalan
series $P$, some coefficients are universal, and, in particular, the
following lemma contains the important value $T_{0,1}=1$.

\Lemma{\label{firstDualCoefficients}
  For any dual coefficients $(T_i)_{i\in\NN^2}$,

  \smallskip

  \noindent (i) $T_{n,0}=0$ for any $n$ in $\NN$;

  \medskip

  \noindent (ii) $\ds T_{n,1}=\cases{1 & if $n=0$,\cr
                                    0 & otherwise.\cr}$
}

\bigskip

The dual coefficients allow us to express the dual basis of a Catalan basis
$(e_k)$ associated to the same Catalan power series $P$. Viewing $P(z,t)$ 
as a power series in $t$, and given that for any nonnegative integer $n$
the ratio $(s^n-t^n)/(s-t)$ is a polynomial in $(s,t)$,
$$
  \Delta P(z,s,t)
  ={P(z,s)-P(z,t)\over s-t}
$$
is a power series in $(z,s,t)$; we extend this definition to
$\partial P(z,s)/\partial s$ on tuples of the form $(z,s,s)$. In
the power series $\Delta P(z,s,t)$ we agree to write the monomial in the
form $z^is^jt^k$, with the variables in this order. If $T$
is an operator then $\Delta P (z,T,t)$ makes sense.

\bigskip

\noindent{\bf Examples.} (continued) a) $\Delta P(z,s,t)=1$.

\noindent b) $\Delta P(z,s,t)=1-z(s+t)$.

\Theorem{\label{dualBasis}
  The dual forms of a Catalan basis $(e_k)_{k\in\NN}$ with
  dual coefficients $(T_i)_{i\in\NN^2}$ are given as follows. Let
  $$
    T(f)(z)=\sum_{i\in\NN^2} T_i \tilde e_i(z,1) 
    f(q^{\langle i,\e_2\rangle}z) \, .
    \eqno{\equa{TExplicit}}
  $$
  Then
  $$
    [e_n]f=q^n[z^0]\Bigl( {\Delta P(z,T,q^n)f(z)\over e_n(qz)P(z,1)}\Bigr)\,.
    \eqno{\equa{dualForms}}
  $$
}

Put differently, this result asserts that if one knows how to decompose
the projection $(z,t)\mapsto t$ over the predual 
basis $(\tilde e_i)_{i\in\NN^2}$, viz.\ \dualCoefficientsEq, then there is 
an explicit formula for the dual forms $[e_n]$.

\bigskip

\noindent{\bf Examples.} (continued) a) We have $Tf(z)=T_{0,1}\tilde e_{0,1}(z,1)
f(qz)=f(qz)$. Thus, 
$$
  [e_n]f=q^n[z^0]{f(z)\over q^n z^n} = [z^n]f(z) \, .
$$ 
\noindent b) We have 
$$\eqalign{
  Tf(z)
  &{}=\sum_{n\in\NN} \calC_n \tilde e_{n,n+1}(z,1)f(q^{n+1}z) \cr
  &{}=\sum_{n\in\NN} \calC_n z^{-1} e_{n+1}(z) f(q^{n+1}z) \, .\cr
  }
$$
Since $\Delta P(z,T,q^n)=1-z(T+q^n)$, we obtain
$$\eqalign{
  [e_n]f
  &{}=q^n [z^0] { f(z)-zTf(z)-zq^nf(z)\over e_n(qz)P(z,1)} \cr
  &{}=[z^0] {f(z)-zTf(z)-zq^nf(z)\over z^n (z,q)_{n+1}} \cr
  &{}=[z^n] {f(z)\over (z,q)_n} -[z^{n-1}]{Tf(z)\over (z,q)_{n+1}} \, . \cr
  }
$$

\bigskip

In general $\Delta P(z,T,q^n)$ involves powers of $T$. Thus, to 
use \dualForms, we need to be able to calculate powers of $T$.
If an operator can be diagonalized, then its powers can be calculated very
efficiently. In the present situation, we will see in section \fixedref{4} that
$T$ is diagonal in the basis $(e_k)$ and more precisely that $Te_k=q^k e_k$. 
Therefore, to use a diagonalization method to calculate $T^nf$ requires one to 
expand $f$ on the basis $e_n$, which is in effect to know the dual basis! A
more effective solution is to use the next result. We define the real-valued
maps $B_n$ on $(\RR^2)^n$ by
$$
  B_n(u_1,\ldots,u_n)
  =\cases{ \sum_{1\leq i<j\leq n} \langle u_i,\e_2\rangle
      \langle u_j,\e_1\rangle & if $n\geq 2$,\cr
      0 & for $n=0,1$.\cr}
  \eqno{\equa{BnDef}}
$$
A simple consideration of $B_3(u_1,u_2,u_3)$ shows that this map is not
multilinear in general, and, comparing with $B_3(u_2,u_1,u_3)$, is not symmetric
in general.

\Proposition{\label{TPower}
  For any integer $n$ positive,
  $$\displaylines{\quad
    T^nf(z)
    =\sum_{i_1,\ldots,i_n\in\NN^2} T_{i_1}\ldots T_{i_n} 
    q^{B_n(i_1,\ldots,i_n)}
    \tilde e_{i_1+\cdots+i_n}(z,1) 
    \hfill\cr\hfill
    f(q^{\langle i_1+\cdots+i_n,\e_2\rangle}z)
    \, .\cr}
  $$
}


\def\prevs{\the\sectionnumber .\the\snumber }
\def\preveq{(\the\sectionnumber .\the\equanumber)}

\section{Some  equations with \poorBold{$q$}-com\-muting coefficients, 
generating functions and recursions}
The purpose of this section is to show that, underpinning Theorem \dualBasis,
there is an interesting connection between some algebraic equations in
some noncommuting variables and some generating functions. In essence,
we show the equivalence between the problem of determining dual bases to 
Catalan bases, determining dual coefficients, solving some
equations in $q$-commuting variables via a $q$-commuting analogue of Puiseux
series, and solving some recusions which generalize the Segner recursion
for Catalan numbers. In short, for a variety of different problems, solving
one allows one to solve all the others.

We will not be interested in general noncommuting variables, but 
in $q$-commuting ones in the following sense.

\Definition{\label{qCommutingVariables}
  Two variables $(A,M)$ $q$-commute if $AM=qMA$.
}

\bigskip
Since we will consider power series of $q$-commuting variables,
we agree on the following in order to avoid any ambiguity.

\bigskip

\noindent{\bf Convention.} {\it 
  (i) When writing power series in a pair of $q$-commuting variable $(X,Y)$, we
  always write the monomials in the form $X^iY^j$.

  \noindent (ii) If $i=(i_1,i_2)$ is in $\NN^2$, we write $(X,Y)^i$ for
  $X^{i_1}Y^{i_2}$.
}

\bigskip

Our next result builds upon Newton's method for finding Puiseux series
expansion of roots of polyomials (see Walker, 1950) and its multivariable
extension of McDonald (1995). It shows that solving some equations in
$q$-commuting coefficients is equivalent to determining dual coefficients.

\Theorem{\label{qEquation}
  Let $(M,A)$ be two $q$-commuting variables. Let $P$ be a Catalan power series
  and let $(T_i)_{i\in\ZZ^2}$ be its dual coefficients. 
  Then
  $$
    T=\sum_{i\in\NN^2} T_i (M,A)^i
  $$
  is the unique power series in $(M,A)$ solving the equation
  $$
    A=P(M,T) \, .
  $$
}


Stated differently, Theorem \qEquation\ asserts that 
if $\sum_{i\in\NN^2} T_i(M,A)^i$ is a solution of $A=P(M,T)$, then
the generating function of the coefficients $(T_i)_{i\in\NN^2}$ with respect to 
the basis $\bigl(\tilde e_i(z,t)\bigr)_{i\in\NN^2}$ is easy to calculate, 
for it is exactly $t$ thanks to \dualCoefficientsEq.

\bigskip

\noindent{\bf Examples.} (continued) a) Since only $T_{0,1}$ does not vanish
and is equal to $1$, Theorem \qEquation\ asserts that $T_{0,1}(M,A)^{(0,1)}=A$
is the only solution of the equation $A=T$.

\noindent b) Assume that $q=1$, and let us switch to lower case symbols to
stress the commutativity. The equation $a=t-mt^2$ has solutions
$$
  t_+={1+\sqrt{1-4am}\over 2m} \qquad\hbox{and}\qquad
  t_-={1-\sqrt{1-4am}\over 2m} \, .
$$
Of those two solutions only $t_-$ is a power series in $(a,m)$ since $t_+$ is
a Laurent series having a term $1/m$ of degree $-1$. Recall that the
Catalan numbers
$$
  C_i={(2i)!\over i!(i+1)!}
$$
have the generating function
$$
  \widehat C(z)=\sum_{i\in\NN} C_i  z^i = {1-\sqrt{1-4z}\over 2z} \, .
$$
Since $t_-=a\widehat C(am)$, we have
$$
  [(a,m)^i]t_-
  =C_{\langle i,\e_2\rangle}\delta_{\langle i,\e_1-\e_2\rangle,1} \, .
$$

When $q$ is not $1$, then Carlitz's $q$-Catalan numbers, $(\calC_n)$, are 
involved. Since by the discussion following \exBb\ the only nonvanishing
dual coefficients are $T_{n,n+1}=\calC_n$, $n\geq 0$,
Theorem \qEquation\ asserts that
$$
  T=\sum_{n\in\NN} \calC_n M^n A^{n+1}
$$
solves $A=T-MT^2$. This provides a new characterization of 
Carlitz's $q$-Catalan numbers in terms of a power series solution of 
a quadratic equation with some $q$-commuting coefficients.

\bigskip

To motivate what follows, recall that the Catalan
numbers obey the Segner recursion
$$
  C_n=\sum_{0\leq i\leq n-1} C_{n-1-i}C_i
  \eqno{\equa{usualSegner}}
$$
with the initial condition $C_0=1$ (see Koshy, 2009). Thus, there is some 
form of recursion
involving the coefficients of $t_-$ in example b). The question arises as to
whether such a recursion exists in the more general context of Theorem 
\qEquation, or, equivalently, for the dual coefficients as defined in
Definition \dualCoefficients. Recall that the maps $B_n$ are defined 
in \BnDef.

\Theorem{\label{Segner}
  Let $P$ be a Catalan power series, and set $R(z,t)=\tilde P(z,t)t$.  
  Let $R_{i,j}=[z^it^j]R(z,t)$.
  The dual coefficients $(T_i)_{i\in\NN^2}$ associated to $P$ 
  are the unique solution to the recursion
  $$\displaylines{
  T_r
  =\One\{\, r=(0,1)\,\} +\sum_{(i,j)\in\NN^2} R_{i,j} 
  \sum_{k_1,\ldots,k_{j+1}\in\NN^2} T_{k_1}\ldots T_{k_{j+1}} 
  \hfill\cr\hfill
  q^{B_{j+1}(k_1,\ldots,k_{j+1})} \One\{\, k_1+\cdots+k_{j+1}=
  r-(i+1,0)\,\} \, .
  \qquad\equa{SegnerRecursion}\cr}
$$
}

If this is not clear that \SegnerRecursion\ is indeed a recursion, the proof
of Theorem \Segner\ makes it clear.

\bigskip

\noindent{\bf Examples.} (continued) a) Since $P(z,t)=t$ we see that $R$ 
is the constant
function $0$. Thus, \SegnerRecursion\ gives $T_r=\One\{\, r=(0,1)\,\}$. Hence,
$T_{0,1}=1$ and all the other $T_{i,j}$ vanish.

\noindent b) Since $R(z,t)=t$, we have 
$$R_{i,j}=\cases{ 1 & if $(i,j)=(0,1)$\cr
                  \noalign{\vskip 2pt}
                  0 & otherwise. \cr}
$$
Recursion \SegnerRecursion\ with $r=(0,m)$ becomes
\hfuzz=2pt
$$\eqalign{
  T_{0,m}
  &{}=\One\{\, m=1\,\} +\sum_{k_1,k_2\in\NN^2} T_{k_1}T_{k_2} q^{B_2(k_1,k_2)}
    \One\{\, k_1+k_2=(-1,m)\,\} \cr
  &{}=\One\{\, m=1\,\} \, . \cr
  }
$$
\hfuzz=0pt
Assume now that we proved that $T_{i,j}=0$ for any $i<n$ and 
any $j\not=i+1$, which
we just did for $n=1$. Then, for $n>1$, \SegnerRecursion\ yields
$$
  T_{n,m}=\sum_{k_1,k_2\in\NN^2} T_{k_1}T_{k_2}q^{B_2(k_1,k_2)}
  \One\{\, k_1+k_2=(n-1,m)\,\} \, .
$$
Using our induction hypothesis
$$\eqalign{
  T_{n,m}
  &{}=\sum_{i,j\in\NN} T_{i,i+1} T_{j,j+1} q^{(i+1)j} 
    \One\{\, i+j=n-1\,;\, i+j+2=m\,\} \cr
  &{}=\sum_{i,j\in\NN} T_{i,i+1} T_{j,j+1} q^{(i+1)j}
    \One\{\, i+j=n-1\,;\, n+1=m\,\} \cr
  }
$$
Thus, $T_{n,m}=0$ if $m\not=n+1$ and
$$
  T_{n,n+1}=\sum_{i,j\in\NN} T_{i,i+1}T_{j,j+1} q^{(i+1)j}
  \One\{\, i+j=n-1\,\} \, .
$$
Writing $\calC_n$ for $T_{n,n+1}$, we obtain
$$
  \calC_n=\One\{\, n=0\,\} +\sum_{i,j\in\NN} \calC_i\calC_j q^{j(i+1)}
  \One\{\, i+j=n-1\,\} \, ,
$$
that is $\calC_0=1$ and for $r\geq 1$,
$$
  \calC_r=\sum_{0\leq i\leq r-1} \calC_i \calC_{r-1-i} q^{(r-1-i)(i+1)} \, .
  \eqno{\equa{qSegner}}
$$
As shown in F\"urlinger and Hofbauer (1985; display preceding (2.2)), this
is the recursion for the Carlitz $q$-Catalan numbers. When
$q$ is $1$, the $\calC_i$ are the usual Catalan numbers and \qSegner\ is
the Segner recursion \usualSegner. Therefore,
recursion \SegnerRecursion\ may be viewed as an abstract form of the Segner 
one.

Note that once we know the result, it is a little easier to calculate the 
dual coefficients from Theorems \qEquation. Indeed, the equation corresponding
to that in Theorem \qEquation\ is $A=T-MT^2$. Knowing that only
the $T_{n,n+1}$ may not vanish, we write $\calC_n$ for $T_{n,n+1}$,
$$
  T=\sum_{i\in\NN} \calC_i M^i A^{i+1} \, .
$$
Substituting this form for $T$ into the equation, we should have
$$\eqalign{
  \sum_{i\in \NN} \calC_i M^i A^{i+1}
  &{}=A+\sum_{i,j\in\NN} \calC_i \calC_j M^{i+1}A^{i+1} M^j A^{j+1}\cr
  &{}=A+\sum_{i,j\in\NN} \calC_i \calC_j q^{j(i+1)} M^{i+j+1} A^{i+j+2} \, . \cr}
$$
Therefore,
$$
  \calC_r
  =\One\{\,  r=0\,\}+\sum_{i,j\in\NN} \calC_i \calC_j q^{j(i+1)}
                     \One\{\, i+j=r-1\,\}
  \, ,
$$
which is \usualSegner\ again.

\bigskip


\section{Proofs} This section contains the proofs of the results stated
in sections \fixedref{2} and \fixedref{3}.

\def\prevs{\the\sectionnumber .\the\subsectionnumber .\the\snumber }
\def\preveq{(\the\sectionnumber .\the\subsectionnumber .\the\equanumber)}

\subsection{Proof of Lemma \ekBasis} Let $P$ be the power series associated 
to $(e_k)$. Since $P$ is a Catalan power series, $P(0,t)=t$.
Therefore, \qCatalanTypeA\ implies $[z^0]\bigl(e_k(qz)/e_k(z)\bigr)=q^k$.
Since $e_k$ is a power series, this implies that it is of order $k$.\hfill\qed

\subsection{Proof of Lemma \representation} Let $(e_k)$ be a 
normalized $q$-Catalan basis associated to $P$. Since $e_k$ is of order $k$
and is normalized, there exists a sequence of power series $(f_k)$ such
that $[z^0]f_k=1$ and $e_k(z)=z^k f_k(z)$. Then
$$
  q^k {f_k(qz)\over f_k(z)}
  = {e_k(qz)\over e_k(z)}
  = {P(z,q^k)\over P(z,1)}
  = q^k {1-z\tilde P(z,q^k)q^k\over 1-z\tilde P(z,1)} \, .
$$
Thus,
$$
  f_k(z)= {1-z\tilde P(z,1)\over 1-z\tilde P(z,q^k)q^k} f_k(qz) \, .
$$
By repeated substitution and using that $f_k(0)=1$, we obtain
$$
  f_k(z)=\prod_{j\geq 0} 
         { 1-q^jz\tilde P(q^jz,1) \over 1-q^{j+k} z\tilde P(q^jz,q^k)} \, .
$$

\noindent (ii) Let $(e_k)$ be a $q$-Catalan basis associated to some
Catalan power series $P$.

If $Q$ is another Catalan power series associated to $(e_k)$ then
$$
  {e_k(qz)\over e_k(z)}
  = {P(z,q^k)\over P(z,1)}
  = {Q(z,q^k)\over Q(z,1)} \, .
$$
Setting $c(z)=Q(z,1)/P(z,1)$, this implies $Q(z,t)=c(z)P(z,t)$.
Since both $P$ and $Q$ are Catalan, their coefficient of $t$ is $1$.
Therefore, $c(z)=1$ and $Q=P$.

\subsection{Proof of Lemma \eTildeBasis} (i) Since $P$ is a Catalan power 
series,
$$
  \prod_{0\leq n<j} P(q^nz,t)
  =t^j \prod_{0\leq n<j}\bigl(1-q^nzt\tilde P (q^nz,t)\bigr)\,.
$$
Therefore, writing $(z,t)^i$ 
for $z^{\langle i,\e_1\rangle}t^{\langle i,\e_2\rangle}$, there exist power 
series $f_i$, $i\in\NN^2$, such that $f_i(0,0)=0$ and
$$
  \tilde e_i(z,t)=(z,t)^i \bigl( 1+f_i(z,t)\bigr) \, .
  \eqno{\equa{eTildeBasisA}}
$$
Since $\tilde e_0=1$ we see that $f_0=0$.
We then need to show that any power series 
$$
  g(z,t)=\sum_{i\in\NN^2} g_i(z,t)^i
  \eqno{\equa{gGeneric}}
$$
can be represented in the basis $\tilde e_i$ as $\sum_{i\in\NN^2}G_i\tilde e_i$.

We now claim that we can determine the $G_i$ recursively by traveling through 
the indices $i$ in $\NN^2$ in the order indicated in the following picture.

\setbox2=\vbox{%
\halign{\hbox to 20pt{\hfill{$\sS#$}}&\hbox to 20pt{\hfill{$\sS#$}}&%
        \hbox to 20pt{\hfill{$\sS#$}}&\hbox to 20pt{\hfill{$\sS#$}}&%
        \hbox to 20pt{\hfill{$\sS#$}}\cr%
\setbox1=\vbox{\hrule width 68pt}\dp1=0pt\wd1=0pt\ht1=0pt%
\raise -2pt\box1 17& 18& 19& 20& 25\cr%
\setbox1=\vbox{\hrule width 48pt}\dp1=0pt\wd1=0pt\ht1=0pt%
\raise -2pt\box1 10& 11& 12& 16& 24\cr%
\setbox1=\vbox{\hrule width 28pt}\dp1=0pt\wd1=0pt\ht1=0pt%
\raise -2pt\box1\phantom{1}%
5&   6&  9& 15& 23\cr%
\setbox1=\vbox{\hrule width 8pt}\dp1=0pt\wd1=0pt\ht1=0pt%
\raise -2pt\box1\phantom{1}%
2&   4&  8& 14& 22\cr%
\setbox1=\hbox{\hskip 25pt\vrule height 18pt depth 0pt%
               \hskip 20pt\vrule height 30pt depth 0pt%
               \hskip 20pt\vrule height 42pt depth 0pt%
               \hskip 19pt\vrule height 54pt depth 0pt}\dp1=0pt\ht1=0pt\wd1=0pt
\box1 1&   3&  7& 13& 21\cr%
}%
}%
\dp2=0pt\ht2=0pt\wd2=0pt%
\setbox3=\hbox{\hskip 142pt\raise 13pt\box2}\dp3=0pt\ht3=0pt\wd3=0pt
\parshape=4 0in 2.6in 0in 2.6in 0in 2.6in 0in 2.6in %
To make\box3\  this explicit, considering the coefficient of $(z,t)^0$, we must 
have $G_0=g_0$.

\parshape=4 0in 2.6in 0in 2.6in 0in 2.6in 0in \hsize
For $i,j$ in $\NN^2$, we write $i\prec j$ if $\langle i,\e_1\rangle$ is at 
most $\langle j,\e_1\rangle$, $\langle i,\e_2\rangle$ is at 
most $\langle j,\e_2\rangle$, and $i$ is not $j$; geometrically, that means 
that $i$ is in the rectangle determined by the origin and the 
point $j$, with the vertex $j$ excluded.

Applying $[(z,t)^j]$ to both sides of \gGeneric,
$$
  g_j=\sum_{i\in\NN^2} G_i [(z,t)^j] \tilde e_i\, .
$$
But \eTildeBasisA\ implies that for $[(z,t)^j]\tilde e_i$ not to be $0$, we
must have $i\prec j$ or $i=j$. Therefore, since \eTildeBasisA\ also
implies that $[(z,t)^j]\tilde e_j=1$, we have
$$
  g_j=\sum_{i\prec j} G_i [(z,t)^j]\tilde e_i + G_j \, .
$$
In other words, we can express $G_j$ in terms of $g_j$ and the $G_i$ with
$i\prec j$. The result is then clear.

\subsection{Proof of Lemma \firstDualCoefficients} Given definitions
\CatalanType, \preDualBasis\ and \dualCoefficients, we have
$$
  \sum_{i\in\NN^2} T_i z^{\langle i,\e_1\rangle} t^{\langle i,\e_2\rangle}
  \prod_{0\leq j<\langle i,\e_2\rangle} \bigl( 1-zt\tilde P(q^jz,t)\bigr)
  = t \, .
  \eqno{\equa{proofFirstDualCoefficientsA}}
$$
We then view both sides of \proofFirstDualCoefficientsA\ as power
series in $t$ with coefficients in $\CC[[\,z\,]]$.

\noindent (i) The coefficient of
$t^0$ is obtained when $\langle i,\e_2\rangle=0$ and is
$\sum_{n\in\NN} T_{n,0}z^n$. Thus, given the right hand side 
of \proofFirstDualCoefficientsA, we obtain $T_{n,0}=0$ for any nonnegative 
integer $n$.

\noindent (ii) The coefficient of $t$ in the left hand side of
\proofFirstDualCoefficientsA\ is obtained when $\langle i,\e_2\rangle=1$
and is
$$
  \sum_{n\in\NN} T_{n,1} z^n [t^0]\bigl( 1-zt\tilde P(z,t)\bigr)
 = \sum_{n\in\NN} T_{n,1}z^n \, .
$$
Given \proofFirstDualCoefficientsA, this series must be $1$ and assertion
(ii) of the lemma follows.\hfill\qed

\subsection{Some equations with \poorBold{$q$}-commuting coefficients and 
recursions} 
In order to prove Theorems \dualBasis,
\qEquation\ and \Segner, we first prove a form of equivalence between
Theorems \qEquation\ and \Segner.

\Theorem{\label{First}
  Let $(A,M)$ be some $q$-commuting variables, and let $P$ be a Catalan power 
  series and set $R(z,t)=\tilde P(z,t)t$. The equation $A=P(M,S)$
  has a unique power series $S=\sum_{i\in\NN^2}S_i(M,A)^i$ solution, which is
  defined by the recursion
  $$\displaylines{\quad
      S_r=\One\{\, r=(0,1)\,\} +\sum_{(i,j)\in\NN^2} R_{i,j}
      \sum_{k_1,\ldots,k_{j+1}\in\NN^2} S_{k_1}\ldots S_{k_{j+1}}
    \hfill\cr\hfill
      q^{B_{j+1}(k_1,\ldots,k_{j+1})} \One\{\, k_1+\cdots+k_{j+1}=r-(i+1,0)\,\}
      \, . \quad\equa{FirstEq}\cr}
  $$
}

\Proof
 The proof consists in considering
a generic power series $S$ in $(M,A)$, and show that the relation $A=P(M,S)$
allows us to calculate recursively the coefficients of the various powers
of $(M,A)$ in $S$. 

The maps $B_n$ defined in \BnDef\ come from the following lemma which 
shows that it
allows us to keep track of the exponent of $q$ when multiplying monomials in 
two $q$-commuting variables.

\Lemma{\label{MAPowers}
  Let $(A,M)$ be some $q$-commuting variables.
  For any $k_1,\ldots ,k_n$ in $\NN^2$,
  $$
    (M,A)^{k_1}\ldots (M,A)^{k_n} 
    = q^{B_n(k_1,\ldots ,k_n)} (M,A)^{k_1+\cdots +k_n} \, .
  $$
}

\Proof The proof is by induction on $n$. The statement is obvious if $n$ 
is $1$. Assume that this identity holds for $n$. Then
$$\displaylines{\qquad
    (M,A)^{k_1}\ldots (M,A)^{k_{n+1}}
    =  q^{B_n(k_1,\ldots,k_n)} M^{\langle k_1+\cdots+k_n,\e_1\rangle} 
  \hfill\cr\hfill
    A^{\langle k_1+\cdots+k_n,\e_2\rangle}
    M^{\langle k_{n+1},\e_1\rangle} A^{\langle k_{n+1},\e_2\rangle} \, .
  \qquad\cr}
$$
Since $(A,M)$ $q$-commute, this is
$$
  q^{B_n(k_1,\ldots ,k_n)+\langle k_{n+1},\e_1\rangle\langle k_1
     +\cdots+k_n,\e_2\rangle} (M,A)^{k_1+\cdots+k_{n+1}} \, .
$$
The exponent of $q$ in this formula is
$$
  B_n(k_1,\ldots,k_n)
  +\sum_{1\leq i<n+1}\langle k_i,\e_2\rangle\langle k_{n+1},\e_1\rangle \, ,
$$
which is indeed $B_{n+1}(k_1,\ldots ,k_{n+1})$.\hfill\qed

\bigskip

We can now prove Theorem \First. We express the equation $A=P(M,S)$ as
$$
  S=A+MR(M,S)S \, ,
  \eqno{\equa{qEquationR}}
$$
meaning
$$
  S=A+\sum_{(i,j)\in \NN^2} R_{i,j} M^{i+1}S^{j+1}\,.
$$
We substitute our tentative expansion for $S$, obtaining
$$\displaylines{\qquad
  \sum_{k\in \NN^2} S_k (M,A)^k
  = A+\sum_{(i,j)\in \NN^2} R_{i,j} \sum_{k_1,\ldots,k_{j+1}\in\NN^2}
    S_{k_1}\ldots S_{k_{j+1}} 
  \hfill\cr\hfill
    M^{i+1} (M,A)^{k_1}\ldots (M,A)^{k_{j+1}} \, .
  \qquad\cr}
$$
Using Lemma \MAPowers, we rewrite this equation as
$$\displaylines{\qquad
  \sum_{k\in \NN^2} S_k (M,A)^k
  = A+\sum_{(i,j)\in \NN^2} R_{i,j} \sum_{k_1,\ldots,k_{j+1}\in\NN^2}
  S_{k_1}\ldots S_{k_{j+1}} 
  \hfill\cr\hfill
  q^{B_{j+1}(k_1,\ldots,k_{j+1})}
  (M,A)^{k_1+\cdots+k_{j+1}+(i+1,0)} \, .\cr}
$$
Equating the coefficients of $(M,A)^r$ on both sides of the identity, we
should have for any $r$ in $\NN^2$,
$$\displaylines{
  S_r
  =\One\{\, r=(0,1)\,\} +\sum_{(i,j)\in\NN^2} R_{i,j} 
  \sum_{k_1,\ldots,k_{j+1}\in\NN^2} S_{k_1}\ldots S_{k_{j+1}} 
  \hfill\cr\hfill
  q^{B_{j+1}(k_1,\ldots,k_{j+1})} \One\{\, k_1+\cdots+k_{j+1}=
  r-(i+1,0)\,\} \, .
  \qquad\equa{recursionS}\cr}
$$
Think of the $S_r$ as placed on the lattice $\NN^2$ with $S_r$ sitting
at site $r$. The above equality asserts that $S_r$ can be calculated
from the various $S_k$ 
with $\langle k,\e_1\rangle\leq\langle r,\e_1\rangle-i-1$,
and so $\langle k,\e_1\rangle <\langle r,\e_1\rangle$.
In other words, we can calculate these coefficients recursively once
we know all the coefficients $S_{n,0}$, $n\geq 0$. To show that all those
vanish, we consider $r$ of the form $(n,0)$. Note that the second indicator
function in \recursionS\ is $0$ if $n-i-1$ is negative. This is the case
in particular if $n=0$. Therefore, \recursionS\ implies $S_{0,0}=0$.

Next, note that if $k_1+\cdots+k_{j+1}=(n-i-1,0)$ and all the $k_i$ are
in $\NN^2$, then necessarily $\langle k_i,\e_2\rangle=0$ for $i=1,\ldots,j+1$;
therefore, for $n$ positive, \recursionS\ yields
$$\displaylines{\qquad
  S_{n,0}
  =\sum_{\scriptstyle 0\leq i\leq n-1\atop\scriptstyle j\in\NN} R_{i,j} 
  \sum_{n_1,\ldots,n_{j+1}\in\NN} S_{n_1,0}\ldots S_{n_{j+1},0} 
  \hfill\cr\hfill
  \One\{\, n_1+\cdots+n_{j+1}= n-i-1\,\} \, .
  \qquad\cr}
$$
This allows us to calculate $S_{n,0}$ from $S_{n-1,0},\ldots,S_{0,0}$, and, if
$S_{n-1,0},\ldots,S_{0,0}$ all vanish, so does $S_{n,0}$. Since we have seen
that $S_{0,0}$ vanishes, all $S_{n,0}$ do as well, and this proves Theorem
\First.

\subsection{A first expression for the dual forms}
The purpose of this subsection is to prove a form of Theorem \dualBasis\
which is of interest for theoretical purposes, notably to prove 
Theorem \dualBasis, but does not give yet a computable form of the dual
basis.

\Theorem{\label{implicitDualBasis}
  Let $(e_k)_{k\in\NN}$ be a Catalan basis with associated series $P$, and let
  $T$ be the linear operator on power series defined by $Te_k=q^ke_k$, 
  $k\in\NN$. Then\note{settle the convergence issues}
  $$
    [e_n]f
    =q^n [z^0]\Bigl( {\Delta P(z,T,q^n)f(z)\over e_n(qz)P(z,1)}\Bigr)
    \, .
    \eqno{\equa{implicitDualForms}}
  $$
}

Note that while this theorem provides the dual basis, it is not explicit. 
Indeed, the operator $T$ is defined on the basis $(e_k)$, and in order
to calculate $Tf$ for arbitrary power series $f$, we need a priori to
decompose $f$ on the basis $(e_k)$, and, at this stage of the proof, we
do not know how to do such a decomposition since we do not know the dual
forms $[e_k]$. However, comparing \implicitDualForms\ with \dualForms, we see
that Theorem \dualBasis\ can be proved by showing that $T$, as defined in
Theorem \implicitDualBasis, can be represented as in \TExplicit.

\bigskip

\Proof Define the operator $\Delta_q$
by
$$
  \Delta_qf(z)=f(qz)-f(z) \, .
$$
This operator has the property that for any power series $f$,
$$
  [z^0]\Delta_qf(z)=0 \, .
$$
Following Hofbauer (1984) and Krattenthaler (1984), 
we consider
$$
  \Delta_q \Bigl({e_k\over e_n}\Bigr)(z)
  ={e_k(z)\over e_n(qz)}
  \Bigl( {e_k(qz)\over e_k(z)}-{e_n(qz)\over e_n(z)}\Bigr) \, .
$$
Therefore, since $(e_k)$ is a Catalan basis with associated series $P$,
$$
  \Delta_q  \Bigl({e_k\over e_n}\Bigr)(z)
  = {e_k(z)\over e_n(qz)P(z,1)} (q^k-q^n) \Delta P(z,q^k,q^n) \, .
$$
This implies that if $k\not=n$ then
$$
  [z^0] \Bigl({e_k(z)\over e_n(qz)P(z,1)} \Delta P(z,q^k,q^n)\Bigr)=0 \, ,
$$
while, if $k=n$,
$$\eqalign{
  [z^0] \Bigl({e_k(z)\over e_n(qz)P(z,1)} \Delta P(z,q^k,q^n)\Bigr)
  &{}= [z^0] \Bigl( {1\over P(z,q^k)} {\partial \over\partial t}P(z,t)_{t=q^n}\Bigr) \cr
  &{}={1\over P(0,q^n)} {\partial\over\partial t} P(0,t)|_{t=q^n}\, .  \cr
  }
$$
Since $P$ is a Catalan power series, 
$$
  P(0,q^n)=q^n \qquad\hbox{and}\qquad {\partial\over\partial t} P(0,t)=1 \, .
$$
Therefore,
$$
  q^n [z^0] 
  \Bigl(  {e_k(z)\over e_n(qz)P(z,1)} \Delta P(z,q^k,q^n)\Bigr)
  = \delta_{k,n} \, .
$$
Since $\Delta P(z,s,t)$ is a power series in $(z,s,t)$,
$$
  e_k(z)\Delta P(z,q^k,q^n) = \Delta P(z,T,q^n)e_k(z)\, .
$$
This is the result.

\subsection{Connecting Theorems \First\ and \implicitDualBasis}
The purpose of this subsection is to show that the operator $T$ involved 
in Theorem \implicitDualBasis\ solves a particular case of the equation
involved in Theorem \First.

\Theorem{\label{TEquation}
  Let $(e_k)$ be a Catalan basis associated to a power series $P$, 
  and define a linear operator $T$ by
  $Te_k=q^ke_k$. Consider the operators
  $$
    Af(z)=P(z,1)f(qz)\qquad\hbox{and}\qquad
    Mf(z)=zf(z) \, .
  $$
  Then,

  \noindent
  (i) $(A,M)$ $q$-commute;

  \noindent
  (ii) $A=P(M,T)$.
}

\bigskip

\Proof Since (i) is trivial, only (ii) needs to be proved.
Since $(e_k)$ is a Catalan basis with associated power series $P$, 
$$
  P(z,1)e_k(qz)=P(z,q^k)e_k(z) \, , 
$$
that is
$$
  A e_k(z)=P(M,T)e_k(z) \, .
$$
Since this equality holds on the basis, the result follows.

\subsection{Proof of Theorem \Segner}
To explain the spirit of the proof, consider an equation $a=P(m,t)$ in
real or complex variable. Assume that there is a power series solution 
$$
  t=\sum_{i\in\NN^2} t_i (m,a)^i \, .
  \eqno{\equa{GFEquation}}
$$
There is a dual viewpoint, which is to consider the array $(t_i)_{i\in\NN^2}$
as a given. It has a generating function 
$\hat t(u,v)=\sum_{i\in\NN^2} t_i (u,v)^i$. Then, this power series satisfies
the equation $a=P\bigl(m,\hat t(m,a)\bigr)$. The proof of Theorem \Segner\ is
a $q$-analogue: starting with the definition of the dual coefficients through
their generating function, our goal is to show that their generating
function satisfies a functional equation which is similar to \GFEquation.

Throughout this subsection, we consider a Catalan power series
with dual coefficients $(T_i)_{i\in\NN^2}$, and we set $R(z,t)=\tilde P(z,t)t$.
Let $(\tilde e_i)_{i\in\NN^2}$ be the corresponding predual basis.
For any positive integer $n$ we consider a map $\Pi_n$, 
which maps $n$ power 
series of two variables to a single power series in two variables, 
is $n$-linear, and is defined by
$$\displaylines{\qquad
  \Pi_n(\tilde e_{i_1},\ldots,\tilde e_{i_n})(z,t)
  \hfill\cr\hfill
  =\tilde e_{i_1}(z,t)
   \tilde e_{i_2}(q^{\langle i_1,\e_2\rangle}z,t)\cdots
   \tilde e_{i_n}(q^{\langle i_1+i_2+\cdots+i_{n-1},\e_2\rangle}z,t) \, ;
  \qquad\cr}
$$
this map is extended by $n$-linearity to power series, agreeing that for
$$
  f_{(j)}(z,t)=\sum_{i\in\NN^2} f_{(j),i} \tilde e_i \, ,
  \qquad 1\leq j\leq n\, ,
$$
we set
$$
  \Pi_n(f_{(1)},\ldots ,f_{(n)})=\sum_{i_1,\ldots,i_n\in\NN^2}
  f_{(1),i_1}\cdots f_{(n),i_n} \Pi_n(\tilde e_{i_1},\ldots,\tilde e_{i_n})
  \, .
$$

The following proposition expresses some important properties of these
$n$-linear maps. If $i$, $j$ are in $\NN^2$, we write $\det(i,j)$ for
the determinant of the matrix whose first column is the column vector $i$ 
and second column  is $j$.

\Proposition{\label{PiNProp}
  (i) $\ds\Pi_n(\tilde e_{i_1},\ldots,\tilde e_{i_n})
      =q^{B_n(i_1,\ldots,i_n)}\tilde e_{i_1+\cdots+i_n}$.

  \medskip

  \noindent$\eqalign{{\it (ii)\ }\ds\Pi_n(f_1,\ldots ,f_n)
  &{}=\Pi_2\bigl( f_1,\Pi_{n-1}(f_2,\ldots,f_{n})\bigr) \cr
  &{}=\Pi_2\bigl( \Pi_{n-1}(f_1,\ldots, f_{n-1}),f_n\bigr) \, .\cr}$

  \smallskip

  \noindent (iii) If $f=\sum_{i\in\NN^2} f_i \tilde e_i$ 
  and $g=\sum_{i\in\NN^2} g_i \tilde e_i$, then 
  $$
    \Pi_2(f,g)
     =\sum_{i,j\in\NN^2} q^{B_2(i,j)+\det(i,j)} g_i f_j \tilde e_{i+j} \, .
  $$
}

Note that assertion (iii) shows that $\Pi_2(f,g)$ is not $\Pi_2(g,f)$
because of the terms $q^{\det(i,j)}$.

\bigskip

\Proof
(i) We proceed by induction. For $n=2$, the definition of $\Pi_2$ and
the definition of the $\tilde e_i$ we have
$$\displaylines{\qquad
  \Pi_2(\tilde e_i,\tilde e_j)(z,t)
  =z^{\langle i,\e_1\rangle}\prod_{0\leq n<\langle i,\e_2\rangle}P(q^nz,t)
  \hfill\cr\hfill
    {}\times(zq^{\langle i,\e_2\rangle})^{\langle j,\e_1\rangle}\prod_{0\leq n<\langle j,\e_2\rangle} P(q^n zq^{\langle i,e_2\rangle},t)\qquad\cr}
$$
The right hand side of this identity is
$$
  q^{\langle i,\e_2\rangle \langle j,\e_1\rangle}
  z^{\langle i+j,\e_1\rangle}\prod_{0\leq n<\langle i+j,\e_2\rangle}
  P(q^nz,t)\, ,
$$
which is indeed $q^{B_2(i,j)}\tilde e_{i+j}(z,t)$. The induction is then 
immediate.

\noindent (ii) It follows by induction from the proof of (i).

\noindent (iii) Using (i), we have
$$
  \Pi_2(f,g)
  = \sum_{i,j\in\NN^2} q^{B_2(i,j)}f_i g_j \tilde e_{i+j}\, .
$$
Writing $i=(i_1,i_2)$ and $j=(j_1,j_2)$, we have
$$\eqalign{
  B_2(i,j)
  &{}=\langle i,\e_2\rangle\langle j,\e_1\rangle\cr
  &{}=B_2(j,i)+i_2j_1-j_2i_1\cr
  &{}= B_2(j,i)-\det(i,j) \,,\cr}
$$
and the result follows.
\hfill\qed

\bigskip

Garsia's (1981) tangled product is defined on power series through the bilinear
mapping
$$
  \Gamma_2(z^i,z^j)=z^i (q^iz)^j = q^{ij}z^{i+j} \, .
$$
The operator $\Pi_2$ defines a tangled product by 
$$
  \Pi_2(\tilde e_i,\tilde e_j)
  =q^{B_2(i,j)} \tilde e_{i+j}\, ,
$$
which is somewhat an analogue of Garsia's, but at
the level of the predual basis. However, while $z^{i+j}$ relates to $z^i$ and
$z^j$ by a simple product, $\tilde e_{i+j}$ does not relate in a simple way
to $\tilde e_i$ and $\tilde e_j$; this precludes us from defining an analogue 
of Garsia's roofing operator in our setting which would trivialize tangled 
products into ordinary ones.

The following result asserts that $\Pi_n(t,\ldots,t)(z,t)=t^n$ but 
expresses this equality in a more readable way.
 
\Proposition{\label{TGenerating}
  The generating function
  $$
    \widehat T(z,t)=\sum_{i\in\NN^2} T_i \tilde e_i(z,t)\, ,
  $$ 
  which by Definition \dualCoefficients\ is $t$, 
  satisfies $\Pi_n(\widehat T,\ldots,\widehat T)=t^n$.
}

\bigskip

\Proof Given \dualCoefficientsEq, we have for any $i$ in $\NN^2$,
$$
  \sum_{j\in\NN^2} T_j \tilde e_j(q^{\langle i,\e_2\rangle}z,t)
  = t \, .
$$
This identity can be multiplied on both sides by $T_i\tilde e_i$ and
then summed over $i$ to obtain
$$
  t^2=\sum_{i_1\in\NN^2} \Bigl(T_{i_1}\tilde e_{i_1}(z,t)\sum_{i_2\in\NN^2}T_{i_2}
  \tilde e_{i_2}(q^{\langle i_1,\e_2\rangle}z,t)\Bigr) \, .
$$
More generally, using the same principle, $\widehat T(z,t)^m$, that is, $t^m$, 
is
\hfuzz=1pt
$$\displaylines{
    \sum_{i_1\in\NN^2}\Biggl( T_{i_1}\tilde e_{i_1}(z,t)\sum_{i_2\in\NN^2}
    \biggl( T_{i_2}\tilde e_{i_2}(q^{\langle i_1,\e_2\rangle}z,t)
    \sum_{i_3\in\NN^2} \Bigl( T_{i_3}\tilde e_{i_3}
    (q^{\langle i_1+i_2,\e_2\rangle}z,t)
  \hfill\cr\hfill
    \ldots\sum_{i_m\in\NN^2} T_{i_m} \tilde e_{i_m}
    (q^{\langle i_1+i_2+\cdots+i_{m-1},\e_2\rangle}z,t)\Bigr)\biggr)
    \ldots\Biggr)
  \quad\cr\noalign{\vskip 3pt}\qquad
   {}=\sum_{i_1,\ldots,i_m\in\NN^2} T_{i_1}T_{i_2}\ldots T_{i_m} \Pi_m(\tilde e_{i_1},\ldots,\tilde e_{i_m})(z,t) 
  \hfill\equa{TGFA}%
  \cr\noalign{\vskip 3pt}\qquad
  {}=\Pi_m(\widehat T,\ldots, \widehat T)(z,t) \, . \hfill\qed\cr
}
$$
\hfuzz=0pt

We can now conclude the proof of Theorem \Segner. Identity \TGFA\ and 
Proposition \PiNProp.(i) yield
$$
  t^m=\sum_{i_1,\ldots,i_m} T_{i_1}\ldots T_{i_m} q^{B_m(i_1,\ldots,i_m)}
  \tilde e_{i_1+\cdots+i_m}(z,t) \, .
  \eqno{\equa{TGFAa}}
$$
Therefore, by uniqueness of the decomposition in the basis $(\tilde e_i)$,
for any $r$ in $\NN^2$,
$$
  \sum_{i_1+\cdots+i_m=r} T_{i_1}\ldots T_{i_m} q^{B_m(i_1,\ldots,i_m)}
  = [\tilde e_r] t^m \, .
$$
Referring to the sums involved in \SegnerRecursion, we then have
$$\displaylines{\quad
  \sum_{(i,j)\in\NN^2} R_{i,j} \sum_{k_1,\cdots,k_{j+1}\in\NN^2} T_{k_1}\ldots
  T_{k_{j+1}} q^{B_{j+1}(k_1,\ldots,k_{j+1})} 
  \hfill\cr\hfill
  \One\{\, k_1+\cdots+k_{j+1}=r-(i+1,0)\,\}
  \qquad\cr\noalign{\vskip 5pt}\hfill 
  {}= \sum_{(i,j)\in\NN^2} R_{i,j} [\tilde e_{r-(i+1,0)}] t^{j+1} \, .
  \qquad\equa{TGFB}\cr}
$$

Since $z^i\tilde e_j(z,t)=\tilde e_{j+(i,0)}(z,t)$, we have 
$$
  [\tilde e_{r-(i+1,0)}]f(z,t)
  =[\tilde e_r] z^{i+1}f(z,t) \, .
$$
Therefore, \TGFB\ is
$$\eqalign{
  \sum_{(i,j)\in\NN^2} R_{i,j} [\tilde e_r](z^{i+1}t^{j+1})
  &{}=[\tilde e_r]\Bigl(\sum_{(i,j)\in\NN^2} R_{i,j} z^{i+1}t^{j+1}\Bigr) \cr
  &{}=[\tilde e_r]\bigl( zR(z,t)t\bigr) \cr
  &{}=[\tilde e_r]\bigl( -P(z,t)+t\bigr) \cr
  &{}=-[\tilde e_r] P(z,t)+[\tilde e_r]t \, . \cr}
$$
Now, \dualCoefficientsEq\ implies $[\tilde e_r]t=T_r$, and since 
$P(z,t)=\tilde e_{0,1}(z,t)$, we also 
have $[\tilde e_r]P(z,t)=\One\{\, r=(0,1)\,\}$. 
Thus, we obtain that
$$\displaylines{\qquad
    \sum_{(i,j)\in\NN^2} R_{i,j} \sum_{k_1,\ldots, k_{j+1}\in\NN^2}
    T_{k_1}\ldots T_{k_{j+1}}
  \hfill\cr\hfill
    q^{B_{j+1}(k_1+\cdots+k_{j+1})}
    \One\{\, k_1+\cdots+k_{j+1}=r-(i+1,0)\,\} 
  \hfill\cr\hfill
    {}=-\One\{\, r=(0,1)\,\} + T_r \, , \qquad\cr}
$$
which is recursion \SegnerRecursion. Since the recursion has a unique solution,
this proves Theorem \Segner.

\subsection{Proof of Theorem \dualBasis}
Theorem \implicitDualBasis\ gives us an expression for the dual basis
in terms of the operator $T$ defined by $Te_k=q^ke_k$. 
Since \implicitDualForms\ and \dualForms\ are the same expression, 
it suffices to show that this operator $T$ acts on function as 
indicated by \TExplicit.

Theorem \TEquation\ asserts that $T$ satisfies the equation $A=P(M,T)$
for the specific operators $A$ and $M$ given in that theorem.

Theorem \First\ asserts that a solution of the equation $A=P(M,T)$ is
given by $\sum_{i\in\NN^2} T_i (M,A)^i$ where the $T_i$ obey the recursion
\SegnerRecursion. Theorem \Segner\ asserts that those $T_i$ are precisely
the dual coefficients. Therefore, it only remains to prove that the
particular solution $\sum_{i\in\NN^2} T_i (M,A)^i$ of the equation
$A=P(M,T)$ is the one we are looking for; or, put differently, that the
specific $T$ we are interested in, which is defined a priori by $Te_k=q^k e_k$,
coincides with the power series $\sum_{i\in\NN^2}T_i(M,A)^i$, and acts on 
functions according to \TExplicit. Thus, we need to prove the following.

\Theorem{\label{final}
  Let $A$ and $M$ be as in Theorem \TEquation, and let $(T_i)$ be
  defined by \dualCoefficientsEq. Furthemore, let
  $T=\sum_{i\in\NN^2} T_i (M,A)^i$. Then 

  \noindent (i) $Te_k=q^ke_k$ for any nonnegative integer $k$;

  \noindent (ii) $\ds Tf(z)=\sum_{i\in\NN^2} T_i \tilde e_i(z,1)
  f(q^{\langle i,\e_2\rangle}z)$.
}

\bigskip

\Proof 
(i) An induction shows that
$$
  A^i f(z)=\prod_{0\leq j<i} P(zq^j,1) f(q^iz) \, .
$$
Since $(e_k)$ is a Catalan basis with associated series $P$,
$$\eqalign{
  e_k(q^iz)
  &{}= {P(q^{i-1}z,q^k)\over P(q^{i-1}z,1)} e_k(q^{i-1}z) \cr
  &{}={\ds \prod_{0\leq j<i} P(q^jz,q^k)\over\ds \prod_{0\leq j<i} P(q^jz,1)}
    e_k(z) \, . \cr 
  }
$$
Consequently,
$$
  A^ie_k(z)=\prod_{0\leq j<i} P(q^jz,q^k) e_k(z) \, ,
$$
and
$$\eqalign{
  \sum_{i\in\NN^2} T_i (M,A)^i e_k(z)
  &{}=\sum_{i\in\NN^2} T_i z^{\langle i,\e_1\rangle} 
    \prod_{0\leq j<\langle i,\e_2\rangle} P(q^jz,q^k)e_k(z) \cr
  &{}=\sum_{i\in\NN^2} T_i \tilde e_i (z,q^k) e_k(z) \, . \cr}
$$
Since $(T_i)_{i\in\NN^2}$ are the dual coefficients,
$\sum_{i\in\NN^2} T_i\tilde e_i(z,q^k)=q^k$.

\noindent (ii) We have
$$\eqalignno{
  Tf(z)
  &{}=\sum_{i\in\NN^2} T_i z^{\langle i,\e_1\rangle} A^{\langle i,\e_2\rangle}
   f(z) \cr
  &{}=\sum_{i\in\NN^2} T_i z^{\langle i,\e_1\rangle} 
    \prod_{0\leq j<\langle i,\e_2\rangle} P(zq^j,1)
    f(q^{\langle i,\e_2\rangle}z) \cr
  &{}=\sum_{i\in\NN^2} T_i \tilde e_i(z,1) f(q^{\langle i,\e_2\rangle}z) \, .
  &\equa{finalA}\cr
  }
$$
This proves Theorem \final.\hfill\qed

\Remark Note that writing $\underline t$ for the power series 
$\underline t(z,t)=t$, that is the projection $(z,t)\mapsto t$, 
\dualCoefficientsEq\ asserts 
that $\underline t=\sum_{i\in\NN^2} T_i\tilde e_i$. Thus, \finalA\ could 
also be rewritten as 
the simpler looking expression $Tf(z)=\Pi_2(\underline t,f)(z,1)$.

\subsection{Proof of Theorem \qEquation} Consider the dual
coefficients $(T_i)_{i\in\NN^2}$. Let $(A,M)$ be, as in Theorem \qEquation,
arbitrary $q$-commuting variables. By Theorem \Segner, $(T_i)$ is the unique
solution of the recursion \SegnerRecursion. Theorem \First\ then gives that
$T=\sum_{i\in\NN^2} T_i (M,A)^i$ is a power series solving $A=P(M,T)$.

Conversely, Theorem \First\ implies that the only power series solution 
in $(M,A)$
of $A=P(M,T)$ is given by $\sum_{i\in\NN^2} S_i (M,A)^i$ where $S_i$ satistifies
\FirstEq. But the recursion \FirstEq\ is \SegnerRecursion\ (just substitute
the letter $T$ for the letter $S$ in \FirstEq). Therefore, by Theorem 
\Segner, $(S_i)_{i\in\NN^2}$ are the dual coefficients, and this proves Theorem
\qEquation.

\subsection{Proof of Proposition \TPower} For $n=1$, this is the definition
of $T$. We then proceed by induction, assuming that the statement is correct 
for $n$. Then $T^{n+1} f(z)$ is
$$\displaylines{\quad
  \sum_{i_2,\ldots,i_{n+1}\in\NN^2} T_{i_2}\cdots T_{i_{n+1}}
  q^{B_n(i_2,\ldots,i_{n+1})} 
  \hfill\cr\hfill
  T\bigl( \tilde e_{i_2+\cdots + i_{n+1}}(z,1)
  f(q^{\langle i_2+\cdots+i_{n+1},\e_2\rangle}z)\bigr) \, .
  \quad\equa{TPowerA}\cr}
$$
But setting $i=i_2+\cdots +i_{n+1}$ and using \TExplicit,
$$\eqalignno{%
   T\bigl(\tilde e_i(z,1) f(q^{\langle i,\e_2\rangle}z)\bigr)
  &{} =\sum_{j\in\NN^2} T_j \tilde e_j(z,1) 
    \tilde e_i(q^{\langle j,\e_2\rangle}z,1)
    f(q^{\langle j+i,\e_2\rangle}z) \quad \cr
  &{} = \sum_{j\in\NN^2} T_j \Pi_2(\tilde e_j,\tilde e_i)(z,1)
    f(q^{\langle j+i,\e_2\rangle}z) \, .\ 
  &\equa{TPowerB}\cr
  }
$$
Using the first assertion of Proposition \PiNProp\ and substituting $i_1$ for
$j$, we obtain that \TPowerB\ is
$$
  \sum_{i_1\in\NN^2} T_{i_1}q^{B_2(i_1,i)}\tilde e_{i_1+i}(z,1)
  f(q^{\langle i_1+i,\e_2\rangle}z) \, .
$$
Substituting this expression in \TPowerA, the result follows from the identity
$$
  B_n(i_2,\ldots,i_{n+1})+B_2(i_1,i_2+\cdots+i_{n+1})
  = B_{n+1}(i_1,\ldots,i_{n+1})\, . \eqno{\qed}
$$

\bigskip


\def\prevs{\the\sectionnumber .\the\snumber }
\def\preveq{(\the\sectionnumber .\the\equanumber)}%

\section{\poorBold{$p$}-ary powers}
In general the dual coefficients form a two dimensional array. However, in 
example b) of section \fixedref{2} and \fixedref{3}, it was in fact
one-dimensional, since only the $T_{n,n+1}$ could not vanish. More can be
said on this situation and this is related to $P(z,t)$ having a specific
form and $e_k$ being similar to Garsia's (1981) powers, whose definition we will
recall.

\Definition{\label{pAryPower}
  Let $p$ be an integer at least $2$. Given a power series $\phi$ of order
  $1$, its $p$-ary powers are the power series
  $$
    \phi_{p,k,q}(z)=z^k\prod_{0\leq j<k(p-1)}\bigl( 1-\phi(q^jz)\bigr) \, ,
    \qquad k\in\NN \, ,
    \eqno{\equa{qAryPowerDef}}
  $$
  with $\phi_{p,0,q}=1$.
}

\bigskip

Garsia (1981) considers powers defined as follows: start with a power series
$\psi$ of order $1$, without loss of generality satisfying $[z]\psi(z)=1$; he 
considers the series 
$$
  \psi_{[k,q]}(z)=\prod_{0\leq j<k} \psi(q^jz) \, .
  \eqno{\equa{GarsiaPowerDef}}
$$
Defining $\phi$ by $\psi(z)=z-z\phi(z)$, we see 
that $\psi_{[k,q]}=q^{k\choose 2}\phi_{2,k,q}$. Thus, up to a normalization,
Garsia's powers coincide with $2$-ary powers. Gessel (1980) deals with the
more general $p$-ary powers (see his section 12) and his functional 
$\alpha(x_1,\ldots,x_n)$ for general noncommuting variables $x_1,\ldots,x_n$
is the analogue to our $B_n(i_1,\ldots,i_n)$ in his setting.

It is possible to interpret $p$-ary powers as Garsia's powers in at least two
different ways. Indeed, if $\phi$ is a power series of order $1$, we can
consider the series
$$
  \psi(z)=z\prod_{0\leq j<p-1} \bigl( 1-\phi(q^jz)\bigr) \, .
$$
For this specific $\psi$,
$$
  \phi_{p,k,q}(z)=q^{-{k\choose 2}(p-1)} \psi_{[k,q^{p-1}]}(z) \, .
$$
However, this interpretation fails if we allow for $q$ to vary, for
$$\displaylines{\quad
  q^{-{k\choose 2}(p-1)} \psi_{[k,r^{p-1}]}(z)
  \hfill\cr\noalign{\vskip 3pt}\hfill
  {}= q^{-{k\choose 2}(p-1)}r^{k\choose 2} z^k \prod_{0\leq i<k}
  \prod_{0\leq j<p-1} \bigl( 1-\phi(q^jr^{i(p-1)}z)\bigr)\quad\cr}
$$
is not $\phi_{p,k,r}(z)$.

A more useful interpretation is obtained as follows. For any integer $p$ 
at least $2$, define the linear operator $K_p$ by 
$$
  K_pf(z)=f(z^{p-1}) \, .
$$
We then have
$$\eqalignno{
  K_p\phi_{p,k,q}(z)
  &{}= z^{k(p-1)} \prod_{0\leq j<k(p-1)} \bigl( 1-\phi(q^jz^{p-1})\bigr) \cr
  &{}= (K_p\phi)_{2,k(p-1),q^{1/(p-1)}}(z) \, .
  &\equa{GarsiaPAry}\cr}
$$
Thus, up to a normalizing factor, $\phi_{p,k,q}(z^{p-1})$ is the Garsia
power $\bigl( \Id(1-K_p\phi)\bigr)_{[k(p-1),q^{1/(p-1)}]}$.
 
Despite these two interpretations of $p$-ary powers as Garsia powers, 
introducing the notion will make our results easier to prove, and we 
will see in the next section that the notion occurs naturally while the
interpretation in terms of Garsia power is somewhat contrived.

Our next result shows that $p$-ary powers are $q$-Catalan bases associated
to some Catalan power series that have a special form, and that the dual
coefficients have a specific sparsity.

\Proposition{\label{pAryDual}
  The $p$-ary powers $\phi_{p,k,q}$ form the normalized $q$-Catalan basis
  associated with the Catalan power series $P(z,t)=t-t\phi(t^{p-1}z)$.
  The dual coefficients are determined by 
  $$
    T_i=0 \hbox{ if }
    i\not\in\bigl\{\, \bigl(n,n(p-1)+1\bigr)\,:\, n\in\NN\,\bigr\}\, ,
  $$
  and
  $$
    \sum_{n\geq 0} q^{-n}T_{n,n(p-1)+1} \phi_{p,n,q}(qz)={1\over 1-\phi(z)} \, .
    \eqno{\equa{CnpDef}}
  $$
}

\Proof Since $\phi$ is of order at least $1$ and $p$ is at least $2$, the power
series $P$ given in the statement is a Catalan one. Since
$$
  {\phi_{p,k,q}(qz)\over\phi_{p,k,q}(z)}
  = q^k {1-\phi(q^{k(p-1)}z)\over 1-\phi(z)}
  = {P(z,q^k)\over P(z,1)} \, ,
$$
the first assertion of the Proposition follows from Lemma \representation.

The predual basis associated to $(\phi_{p,k,q})$ is then
$$
  \tilde e_i(z,t)
  =z^{\langle i,\e_1\rangle}t^{\langle i,\e_2\rangle}\prod_{0\leq j<\langle i,\e_2\rangle}
    \bigl(1-\phi(q^jt^{p-1}z)\bigr) \, .
$$
Therefore,
$$
  \tilde e_i(z^{p-1},t)
  =z^{(p-1)\langle i,\e_1\rangle} t^{\langle i,\e_2\rangle}
  \prod_{0\leq j<\langle i,\e_2\rangle} 
  \Bigl( 1-\phi\bigl(q^j(tz)^{p-1}\bigr)\Bigr) \, .
$$
By definition of the dual coefficients, $t=\sum_{i\in\NN^2} T_i\tilde e_i(z,t)$;
substituting $z^{p-1}$ for $z$ in this identity and then multiplying both sides
by $z$,
$$\eqalignno{
  zt
  &{}=\sum_{i\in\NN^2} T_i z\tilde e_i(z^{p-1},t) \cr
  &{}=\sum_{i\in\NN^2} T_i 
    z^{(p-1)\langle i,\e_1\rangle+1-\langle i,\e_2\rangle}
    (zt)^{\langle i,\e_2\rangle}
  \prod_{0\leq j<\langle i,\e_2\rangle}
    \Bigl( 1-\phi\bigl(q^j(tz)^{p-1}\bigr)\Bigr) \, .
  \cr
  &&\equa{pAryDualA}\cr}
$$
Setting $\tau=zt$ we see from \pAryDualA\ that whenever $T_i$ does not vanish
we must have
$$
  (p-1)\langle i,\e_1\rangle+1-\langle i,\e_2\rangle =0 \,. 
$$
Thus, only $T_{n,n(p-1)+1}$, $n\in\NN$, may not vanish. Consequently, 
\pAryDualA\ yields
$$
  \tau=\sum_{n\in\NN} T_{n,n(p-1)+1} \tau^{n(p-1)+1}
  \prod_{0\leq j<n(p-1)+1} \bigl( 1-\phi(q^j\tau^{p-1})\bigr) \, ,
$$
that is, after simplifying both sides of this identity by $\tau$ and
setting $z=\tau^{p-1}$,
$$\eqalignno{
  1
  &=\sum_{n\in\NN} T_{n,n(p-1)+1} z^n \prod_{0\leq j<n(p-1)+1} 
    \bigl( 1-\phi(q^jz)\bigr) \cr
  &{}=\bigl( 1-\phi(z)\bigr) \sum_{n\in\NN} q^{-n} T_{n,n(p-1)+1}
    \phi_{p,n,q}(qz)\, .
  &\qed\cr}
$$

Starting with a power series $\phi(z)$ of order $1$, we consider the
$p$-ary powers $\phi_{p,k,q}(z)$; they induce a linear operator $U_{p,\phi,q}$
defined by
$$
  U_{p,\phi,q}z^k = q^{{k(p-1)\choose 2}/(p-1)} \phi_{p,k,q}(z) \, .
$$
Following Garsia (1981) and Krattenthaler (1988), extending their result
from Garsia powers to $p$-ary powers, we will
now construct the inverse of $U_{p,\phi,q}$. 

Given Proposition \pAryDual\ the dual coefficients for the 
basis $(\phi_{p,n,q})$ are defined by \CnpDef.
We define
$$
  \phi^{\diamond p}(z)
  =-\sum_{n\geq 1} T_{n,n(p-1)+1} q^{-{n(p-1)+1\choose 2}/(p-1)} z^n \, .
$$
Note that the map $\phi\mapsto \phi^{\diamond p}$ depends on $q$, a dependence
lost in the notation $\phi^{\diamond p}$.
In order to avoid any ambiguity in the notation, we agree that the mapping
$\phi\mapsto\phi^{\diamond p}$ has the precedence over $\phi\mapsto \phi_{p,k,q}$,
so that, for instance, $\phi^{\diamond p}_{p,k,q}$ 
means $(\phi^{\diamond p})_{p,k,q}$. 

\Theorem{\label{inverseSequence} 
  $U_{p,\phi,q}$ and $U_{p,\phi^{\diamond p},1/q}$ are
  inverse of each others.
}

\bigskip

\Proof We will first prove the result when $p$ is $2$, which corresponds
to the Garsia powers, and extend this result to arbitrary $p$ via \GarsiaPAry.

\noindent{\it Case $p=2$.}  We assume for the time being that $p$ is $2$.
The heart of the proof is identity \TGFAa\ and the following lemma.

\Lemma{\label{BnExpand} 
  Let $(u_i)_{1\leq i\leq m}$ be $m$ vectors in $\NN^2$ with
  $$
    \langle u_i,\e_2\rangle =\langle u_i,\e_1\rangle+1
  $$
  for $i=1,\ldots,m$. Set $v=u_1+\cdots+u_m$. Then
  $$\displaylines{\quad
    B_m(u_1,\ldots,u_m)
    \hfill\cr\noalign{\vskip 3pt}\hfill
    {}= {\langle v,\e_2\rangle\choose 2} 
     -\sum_{1\leq i\leq m} {\langle u_i,\e_2\rangle\choose 2}
     -\sum_{1\leq i\leq m} \langle u_i,\e_2\rangle (m-i) \, .
    \quad\cr}
  $$
}

\Proof Write $n_i$ for $\langle u_i,\e_2\rangle$. Since by hypothesis 
$\langle u_j,\e_1\rangle=n_j-1$,
$$\eqalign{
  B_m(u_1,\ldots,u_m)
  &{}=\sum_{1\leq i<j\leq m} n_i (n_j-1)  \cr
  &{}={1\over 2} \Bigl( \sum_{1\leq i\leq m} n_i\Bigr)^2
    -{1\over 2} \sum_{1\leq i\leq m} n_i^2
    - \sum_{1\leq i\leq m}  n_i (m-i) \cr
  &{}={1\over 2} \langle v,\e_2\rangle^2-{1\over 2} \sum_{1\leq i\leq m} 
    n_i^2-\sum_{1\leq i\leq m} n_i (m-i) \, .\cr}
$$
Note that for any integer $n$,
$$
  {n^2\over 2} ={n\choose 2} +{n\over 2} \, .
$$
Therefore, $B_m(u_1,\ldots,u_m)$ is
$$
  {\langle v,\e_2\rangle\choose 2} + {\langle v,\e_2\rangle\over 2}
  -\sum_{1\leq i\leq m} {n_i\choose 2}
  -{1\over 2} \sum_{1\leq i\leq m} n_i 
  -\sum_{1\leq i\leq m} n_i (m-i) \, .
$$
This is the result since 
$\langle v,\e_2\rangle=\sum_{1\leq i\leq m} n_i$.\hfill\qed

\bigskip

Throughout the proof we will use the $q$-Catalan basis $(e_k)=(\phi_{2,k,q})$,
and switch freely between the notation $e_k$ and $\phi_{2,k,q}$ according
to the context. In particular $[e_k]$ is the dual form of $\phi_{2,k,q}$.
The predual basis associated to $(e_k)$ is
$$
  \tilde e_i(z,t)
  =z^{\langle i,\e_1-\e_2\rangle} e_{\langle i,\e_2\rangle}(zt) \, .
$$
Proposition \pAryDual\ shows that the dual coefficients are determined by
$$
  \sum_{n\geq 0} q^{-n} T_{n,n+1} \phi_{2,n,q}(qz)={1\over 1-\phi(z)} \, .
$$
Identity \TGFAa\ is then
$$\displaylines{\qquad
  t^m
  = \sum_{i_1,\ldots,i_m\in\NN^2} T_{i_1}\ldots T_{i_m}
  q^{B_m(i_1,\ldots,i_m)} z^{\langle i_1+\cdots+i_m,\e_1-\e_2\rangle}
  \hfill\cr\hfill
  e_{\langle i_1+\cdots+i_m,\e_2\rangle} (zt) \, .
  \qquad\cr}
$$
Multiplying both sides by $z^m$ and substituting $\zeta$ for $zt$,
$$\displaylines{\qquad
  \zeta^m
  =\sum_{i_1,\ldots,i_m\in\NN^2} T_{i_1}\ldots T_{i_m} q^{B_m(i_1,\ldots,i_m)}
  z^{\langle i_1+\cdots+i_m,\e_1-\e_2\rangle+m}
  \hfill\cr\hfill
  e_{\langle i_1+\cdots+i_m,\e_2\rangle}(\zeta) \, . 
  \qquad\equa{inverseSequenceA}\cr}
$$
Given Proposition \pAryDual\ with $p=2$, only the $T_i$ for 
which $i=(n-1,n)$, $n\geq 1$, 
may not vanish. Moreover, writing $i_j$ as $(n_j-1,n_j)$, we see that
$$
  \langle i_1+\cdots+i_m,\e_1-\e_2\rangle + m = 0 \, ,
$$
and Lemma \BnExpand\ gives
$$
  B_m(i_1,\ldots,i_m)
  ={n_1+\cdots+n_m\choose 2} 
   -\sum_{1\leq i\leq m} {n_i\choose 2}
   -\sum_{1\leq i\leq m} n_i(m-i) \, . 
$$
Thus, \inverseSequenceA\ gives
$$\displaylines{\quad
  \zeta^m=\sum_{n_1,\ldots,n_m\geq 1} T_{n_1-1,n_1}\ldots T_{n_m-1,n_m}
  \hfill\cr\hfill
  q^{{n_1+\cdots+n_m\choose 2}-\sum_{1\leq i\leq m} {n_i\choose 2}
     -\sum_{1\leq i\leq m} n_i(m-i)} e_{n_1+\cdots+n_m}(\zeta) \, . 
  \cr}
$$
Substituting $z$ for $\zeta$ and applying $[e_k]$ to both sides,
$$\displaylines{\qquad
  [e_k]z^m
  = \sum_{n_1,\ldots,n_m\geq 1} \One\{\, n_1+\cdots+n_m=k\,\} 
  T_{n_1-1,n_1}\ldots T_{n_m-1,n_m}
  \hfill\cr\hfill
  q^{{k\choose 2} -\sum_{1\leq i\leq m} {n_i\choose 2} -\sum_{1\leq i\leq m} n_i(m-i)} 
  \, . \cr}
$$
The key of the proof is to multiply both sides of this identity by $z^k$ and
notice the nice way to factor the right hand side, obtaining
$$\displaylines{\quad
  z^k [e_k]z^m
  = q^{k\choose 2} \sum_{n_1,\ldots,n_m\geq 1} \One\{\, n_1+\cdots+n_m=k\,\}
  \hfill\cr\noalign{\vskip -7pt}\hfill
  \prod_{1\leq i\leq m} T_{n_i-1,n_i} q^{-{n_i\choose 2}} 
  \Bigl({z\over q^{m-i}}\Bigr)^{n_i} \, .
  \cr}
$$
We then multiply both sides by $q^{-{k\choose 2}}$ and sum over $k$ and use
$T_{0,1}=1$ to obtain, with $\theta(z)=z\bigl( 1-\phi^{\diamond 2}(z)\bigr)$, 
that
$$\eqalignno{
  \sum_{k\geq 0} q^{-{k\choose 2}} z^k [e_k]z^m
  &{}=\prod_{1\leq i\leq m} \biggl(\sum_{n\geq 1} T_{n-1,n} q^{-{n\choose 2}} 
    \Bigl({z\over q^{m-i}}\Bigr)^n\biggr) \cr
  &{}=\prod_{1\leq i\leq m} \theta \Bigl( {z\over q^{m-i}}\Bigr) \cr
  &{}=U_{2,\phi^{\diamond 2},1/q} z^m \, .
  &\equa{inverseSequenceB}\cr}
$$
We apply $U_{2,\phi,q}$ on both sides of the identity and use that $([e_k])$ 
is the dual basis of $(\phi_{2,k,q})$ to get
$$
  z^m=U_{2,\phi,q}U_{2,\phi^{\diamond 2},1/q} z^m \, ,
$$
and therefore $U_{2,\phi,q}U_{2,\phi^{\diamond 2},1/q}=\Id$.

To prove that $U_{2,\phi^{\diamond 2},1/q}U_{2,\phi,q}=\Id$, we can either use 
that since both $U_{2,\phi,q}$ and $U_{2,\phi^{\diamond 2},1/q}$ map basis to basis
they are invertible and therefore their right and left inverse coincide; or,
we can note that \inverseSequenceB\ asserts that
$$
  \sum_{k\geq 0} q^{-{k\choose 2}} z^k [e_k](z^\ell)
  = q^{-{\ell\choose 2}}\phi^{\diamond 2}_{2,\ell,1/q}(z) \, .
  \eqno{\equa{inverseSequenceC}}
$$
We then have
$$\eqalignno{
  U_{2,\phi^{\diamond 2},1/q}U_{2,\phi,q}z^m
  &{}=U_{2,\phi^{\diamond 2},1/q} q^{m\choose 2}
    \sum_{\ell\in\NN} [z^\ell]( \phi_{2,m,q}) z^\ell \cr
  &{}=q^{m\choose 2}\sum_{\ell\in\NN} [z^\ell] (\phi_{2,m,q}) 
    q^{-{\ell\choose 2}}\phi^{\diamond 2}_{2,\ell,1/q}(z) \,. \qquad
  &\equa{inverseSequenceD}\cr}
$$
Thus, using \inverseSequenceC, the right hand side of \inverseSequenceD\ is
$$
    q^{m\choose 2}\sum_{k,\ell\in\NN} [z^\ell] (\phi_{2,m,q}) 
    q^{-{k\choose 2}} z^k [e_k](z^\ell)
    =q^{m\choose 2}\sum_{k\in\NN} q^{-{k\choose 2}} 
    z^k [e_k] \phi_{2,m,q} \, . 
$$
Since $\phi_{2,m,q}=e_m$, we obtain that the right hand side 
of \inverseSequenceD\ is
$$
  q^{m\choose 2}\sum_{k\in\NN} q^{-{k\choose 2}} z^k \delta_{k,m}
  = z^m \, , 
$$
and therefore, $U_{2,\phi^{\diamond 2},1/q}U_{2,\phi,q}z^m=z^m$. This proves Theorem 
\inverseSequence\ when $p$ is $2$.

\noindent{\it Case $p\geq 3$}. We assume now that $p$ is at least $3$. 
We set $r=q^{1/(p-1)}$. The identity \GarsiaPAry\ asserts that
$$
  K_p U_{p,\phi,q}z^k = U_{2,K_p\phi,r}K_pz^k \, .
$$
Given its left hand side, this equality shows that the image 
of $U_{2,K_p\phi,r}K_p$ is in the range of $K_p$. Therefore,
$$
  U_{p,\phi,q}=K_p^{-1} U_{2,K_p\phi,r}K_p \, ,
$$
and, consequently, using the present theorem with $p=2$,
$$
  U_{p,\phi,q}^{-1}
  = K_p^{-1} U_{2,K_p\phi,r}^{-1}K_p
  = K_p^{-1} U_{2,(K_p\phi)^{\diamond 2},1/r}K_p \, .
  \eqno{\equa{inverseSequenceF}}
$$
For this formula to be useful, we need to calculate $(K_p\phi)^{\diamond 2}$.
Thus, we consider the $2$-ary powers $(K_p\phi)_{2,k,r}$, and, in order to avoid
any confusion, we write $(S_i)_{i\in\NN^2}$ for their dual coefficients. 
Applying Proposition \pAryDual, these dual coefficients are defined by
$$
  \sum_{n\in\NN} r^{-n} S_{n,n+1} (K_p\phi)_{2,n,r}(rz)
  ={1\over 1-K_p\phi(z)}\, .
  \eqno{\equa{inverseSequenceE}}
$$
Since $(K_p\phi)_{2,n,r}$ is in $z^n\CC[[z^{p-1}]]$ 
and $1/\bigl( 1-K_p\phi(z)\bigr)$ is in $\CC[[z^{p-1}]]$, \inverseSequenceE\
implies that the coefficients $S_{n,n+1}$ must vanish if $n$ is not a multiple
of $p-1$. Therefore, only $S_{n(p-1),n(p-1)+1}$, $n\in\NN$, may not vanish.
It follows that
$$
  (K_p\phi)^{\diamond 2}(z)
  =-\sum_{n\geq 1} S_{n(p-1),n(p-1)+1}r^{-{n(p-1)+1\choose 2}} z^{n(p-1)} \, .
  \eqno{\equa{inverseSequenceG}}
$$
Moreover, since only $S_{n(p-1),n(p-1)+1}$, $n\in\NN$, may not vanish,
\inverseSequenceE\ can then be rewritten as
$$\displaylines{\quad
  \sum_{n\in\NN} S_{n(p-1),n(p-1)+1} z^{n(p-1)}
  \prod_{0\leq j<n(p-1)} \Bigl( 1-\phi\bigl( (zr^{j+1})^{p-1}\bigr)\Bigr)
  \hfill\cr\hfill
  {}= {1\over 1-\phi(z^{p-1})} \, .\cr}
$$
Substituting $z$ for $z^{p-1}$, this means
$$
  \sum_{n\in\NN} q^{-n} S_{n(p-1),n(p-1)+1} \phi_{p,n,q}(qz)
  = {1\over 1-\phi(z)} \, .
$$
Consequently, given \CnpDef, $S_{n(p-1),n(p-1)+1}=T_{n,n(p-1)+1}$ and 
\inverseSequenceG\ becomes
$$
  (K_p\phi)^{\diamond 2} (z)
  = -\sum_{n\geq 1} T_{n,n(p-1)+1} r^{-{n(p-1)+1\choose 2}} z^{n(p-1)} \, .
$$
Given how $\phi^{\diamond p}$ is defined, this means 
$(K_p\phi)^{\diamond 2}(z)=\bigl(K_p(\phi^{\diamond p})\bigr)(z)$. Going back
to \inverseSequenceF, this implies
$$\displaylines{\quad
  K_p^{-1}U_{2,(K_p\phi)^{\diamond 2},1/r} K_p z^m
  \hfill\cr\noalign{\vskip 5pt}\hfill
  \eqalign{
  {}={}& K_p^{-1} r^{-{m(p-1)\choose 2}} 
         \bigl(K_p(\phi^{\diamond p})\bigr)_{2,m(p-1),1/r}(z) \cr
  {}={}&r^{-{m(p-1)\choose 2}} K_p^{-1} 
        \Biggl( z^{m(p-1)}\prod_{0\leq j<m(p-1)}
                \biggl( 1-\phi^{\diamond p}
                       \Bigl({z^{p-1}\over r^{(p-1)j}}\Bigr)
                \biggr)
        \Biggr)\cr
  {}={}&r^{-{m(p-1)\choose 2}} z^m \prod_{0\leq j<m(p-1)} 
        \Bigl( 1-\phi^{\diamond p}\Bigl({z\over q^j}\Bigr)\Bigr) \cr
  {}={}&U_{p,\phi^{\diamond p},1/q} z^m \, .\cr
  }\cr}
$$
This proves Theorem \inverseSequence\ by \inverseSequenceF.\hfill\qed

\bigskip

\Corollary{\label{involution}
  For every $p$ at least $2$, the map $(\phi,q)\mapsto(\phi^{\diamond p},1/q)$ 
  is an involution.
}

\bigskip

\Proof Theorem \inverseSequence\ implies that 
$U_{p,\phi,q}^{-1}=U_{p,\phi^{\diamond p},1/q}$.
Substituting $(\phi^{\diamond p},1/q)$ for $(\phi,q)$, we 
obtain $U_{p,\phi^{\diamond p},1/q}^{-1}
=U_{p,(\phi^{\diamond p})^{\diamond p},q}$. Therefore, 
$U_{p,(\phi^{\diamond p})^{\diamond p},q}=U_{p,\phi,q}$,
which forces $(\phi^{\diamond p})^{\diamond p}=\phi$.\hfill\qed

\bigskip

Starting with a power series $\phi(z)$ of order $1$ and $[z]\phi(z)=1$, 
the power series $P(z,t)=\phi(zt)/z$ is a Catalan one. For this specific 
type of Catalan power series, 
define $P^{\diamond p}(z,t)=\phi^{\diamond p}(zt)/z$. Corollary \involution\ 
implies that $(P,q)\mapsto (P^{\diamond p},1/q)$ is an involution as well.

\bigskip

To make explicit the connection between Theorem \inverseSequence\ and
the work of Garsia (1981) and Krattenthaler (1988), let $\phi$ be a 
power series of order $1$ and set $\psi(z)=z\bigl( 1-\phi(z)\bigr)$. 
Krattenthaler's (1988) identity (6.11) asserts that the relation
$$
  \sum_{k\in\NN} a_k z^k = \sum_{k\in\NN} b_k \prod_{0\leq j<k} \psi(q^j z) \, ,
$$
that is, 
$$
  \sum_{k\in\NN} a_kz^k=U_{2,\phi,q}\sum_{k\in\NN} b_k z^k \, ,
$$ 
is equivalent to
$$
  \sum_{k\in\NN} b_k z^k=U_{2,\phi^{\diamond 2},1/q}\sum_{k\in\NN} a_kz^k\, ,
$$
which is clear from Theorem \inverseSequence.

Similarly, Krattenthaler's (1988) identity (6.10) states that if $F_{k,\ell}$
are defined by
$$
  z^\ell
  =\sum_{k\geq \ell} F_{k,\ell} \prod_{0\leq j<k} \psi(q^jz) 
  =U_{2,\phi,q}\sum_{k\geq \ell} F_{k,\ell} z^k \, ,
$$
then, setting $\Psi(z)=z\bigl( 1-\phi^{\diamond 2}(z)\bigr)$,
$$
  \sum_{k\geq \ell} F_{k,\ell} z^k
  = U_{2,\phi,q}^{-1} z^\ell
  = U_{2,\phi^{\diamond 2},1/q} z^\ell
  = \prod_{0\leq j<\ell} \Psi(z/q^j) \, .
$$

\bigskip

Theorems \qEquation\ and \Segner\ can be viewed as different definitions
of the dual coefficients, in which case Definition \dualCoefficients\ asserts 
that a form of generating function of the dual coefficients, 
$\sum_{i\in\NN^2} T_i \tilde e_i(z,t)$, has a particular value. However, 
this 'form' of generating
function is not as useful as a proper generating function. While we do not 
know how to calculate the generating function of the dual coefficients in
general, our next result shows that in the particular context of $p$-ary
powers, the generating function of the suitably multiplied dual coefficients
obeys a $q$-functional relation, and, as this will be clear in the next 
paragraph, this is related to Theorem \inverseSequence\ and its proof.

Given a power series $\phi$ of order $1$ and the corresponding $p$-ary
powers $\phi_{p,k,q}$, Proposition \pAryDual\ asserts that only 
the $T_{n,n(p-1)+1}$, $n\in\NN$, may not vanish. The generating function of
these dual coefficients is then $\sum_{n\in\NN} T_{n,n(p-1)+1} z^n$. Instead,
it will be easier to consider the generating function 
of $q^{-(p-1){n\choose 2}}T_{n,n(p-1)+1}$,
$$
  \calT_{\phi,p}(z)=\sum_{n\in\NN} q^{-(p-1){n\choose 2}} T_{n,n(p-1)+1}z^n \, .
$$
As an indication that  this generating function is related to 
Theorem \inverseSequence, note the identity
$$\eqalign{
  \calT_{\phi,p}(z)
  &{}=1+\sum_{n\geq 1} q^{-{n(p-1)+1\choose 2}/(p-1)} q^{np/2}
     T_{n,n(p-1)+1}z^n\cr
  &{}=1-\phi^{\diamond p}(q^{p/2}z) \, .\cr
  }
$$

For the following result, it is of interest to use Garsia's roofing operator,
defined on power series $f(z)=\sum_{n\in\NN} f_n z^n$ by
$$
  \widehat f(z)=\sum_{n\in\NN} q^{-{n\choose 2}} f_n z^n \, ,
$$
as well as his reciprocal staring operator
$$
  \lstar f(z)=\prod_{n\in\NN} f(z/q^n) \, .
$$
The roofing operator may be itereated, setting $f\,\widehat{\ }^{\,p}
=\widehat{f\,\widehat{\ }^{\,(p-1)}}$, or, more explicitly,
$$
  f\,\widehat{\ }^{\,p}(z)
  = \sum_{n\in\NN} q^{-p{n\choose 2}} f_n z^n \, .
$$

\Theorem{\label{dualGF}
  Let $\phi(z)=\sum_{i\geq 1} \phi_i z^i$ be a power series of order $1$. 
  The generating function $\calT_{\phi,p}$ obeys the functional relation
  $$
    \calT_{\phi,p}(z) =
    1 +\sum_{i\geq 1} \phi_i z^i q^{-(p-1){i\choose 2}} 
    \prod_{0\leq j\leq (p-1)i}\calT_{\phi,p}(q^{-j}z) \, .
    \eqno{\equa{dualGFEquation}}
  $$
  and satisfies the identity
  $$
    \calT_{\phi,p}(z)
    ={ \bigl(\lstar \,(1-\phi)\big)\widehat{\ }^{\,(p-1)}(z/q)
       \over
       \bigl(\lstar \,(1-\phi)\big)\widehat{\ }^{\,(p-1)}(z) } \, .
  $$
}

Before proving this result, some remarks and an example may help to 
understand its content.

\Remark Defining
$$
  \Phi_q(z)=\prod_{n\in\NN} \bigl(1-\phi(q^nz)\bigr) \, ,
$$
we see that ${}^*(1-\phi)(z)=\Phi_{1/q}(z)$. Set
$$
  g(z)=\sum_{n\in\NN} q^{-(p-1){n\choose 2}} z^n [z^n]\Phi_{1/q}(z) \, ,
$$
that is $g(z)=\Phi_{1/q}\widehat{\ }^{(p-1)}$. The second identity in 
Theorem \dualGF\ gives that $\calT_{\phi,p}(z)= g(z/q)/g(z)$.

\bigskip

\noindent{\bf Example.} Continuing example b) of sections 2 and 3, and 
given Proposition \pAryDual, the Catalan power series $P(z,t)=t-t^2z$ 
corresponds to $\phi(z)=z$ and $p=2$. In this case, $[z]\phi(z)=1$, that is,
$\phi_1=1$, and all the other $\phi_i$ vanish. The assertion of 
Theorem \dualGF\ is then
$$
  \calT_{\phi,2}(z)=1+z\calT_{\phi,2}(z/q)\calT_{\phi,2}(z) \, .
$$
This implies
$$
  \calT_{\phi,2}(z)={1\over 1-z\calT_{\phi,2}(z/q)} \, ,
$$
and, by iterated substitution, we obtain the continued fraction,
$$
  \calT_{\phi,2}(z)
  ={1\over \ds 1-
     {\ds z/q^0\over\ds 1-
       {\ds z/q^1 \over\ds 1-
         {\ds z/q^2\over\ds 1-
           {\ds z/q^3\over\ds \ldots}}}}} \, ,
  \eqno{\equa{continuousFractionT}}
$$
as indicated in Garsia (1981). 

With the notation of the remark preceding this example, we have here
$$
  \Phi_{1/q}(z)=\prod_{n\in\NN} (1-z/q^n) = (z,1/q)_\infty \, .
$$
Using Euler's identity (see Andrews, Askey and Roy, 1999, Corollary 10.2.2),
$$
  [z^n]\Phi_{1/q}(z)=(-1)^n {q^{-{n\choose 2}}\over (1/q,1/q)_n}
$$
Thus, since $p=2$,
$$
  g(z)=\sum_{n\in\NN} (-1)^n {q^{-n(n-1)}\over (1/q,1/q)_n} z^n \, .
$$
Recall that, following Ismail (2009, formula (21.7.3)) the $q$-Airy 
function is defined by
$$
  \Ai_q(z)=\sum_{n\geq 0} (-1)^n {q^{n^2}z^n\over (q,q)_n} \, .
$$
Thus $g(z)=\Ai_{1/q}(qz)$. Theorem \dualGF\ and the remark following it imply 
that in this case
$$
  \calT_{\phi,2}(z)={\Ai_{1/q}(z)\over \Ai_{1/q}(z/q)} \, ,
  \eqno{\equa{AiryT}}
$$
the identity between \continuousFractionT\ and \AiryT\  being the 
Rogers-Ramanujan continuous fraction.

\bigskip

\Remark It is interesting to compare Theorem \dualGF\ with Gessel's (1980)
Theorem 12.2. When $q$ is $1$, \dualGFEquation\ is
$$
  \calT_{\phi,p}(z)=1+\calT_{\phi,p}(z)\phi\bigl(z\calT_{\phi,p}(z)^{p-1}\bigr)
  \, .
  \eqno{\equa{dualGFRemA}}
$$
The function $\calT_{\phi,p}$ can then be obtained by Lagrange inversion.
Setting $y=z\calT_{\phi,p}(z^{p-1})$ and $g(y)=1/\bigl(1-\phi(y)\bigr)$, 
substituting $z^{p-1}$ for $z$ in \dualGFRemA, we
can rewrite \dualGFRemA\ as $y=zg(y^{p-1})$. This is to compare with
(12.1) in Gessel (1980), which, when $q$ is $1$, would be $y=g(zy^{p-1})$. 
The second equality in Theorem \dualGF\ can then be interpreted as a 
Lagrange type inversion formula for the functional equation \dualGFEquation.

To further compare Theorem \dualGF\ with Gessel's (1980) Theorem 12.2, 
setting $r=1/q$ and $m=p-1$, and writing now $f$ for $\calT_{\phi,p}$ and 
recalling the notation in \GarsiaPowerDef, \dualGFEquation\ is
$$\eqalignno{
  f(z)
  &{}= 1+f(z)\sum_{i\geq 1} \phi_i z^i r^{(p-1){i\choose 2}}f_{[mi,r]}(rz)\cr
  &{}=1+f(z)\sum_{i\geq 1} r^{-mi/2}\phi_i r^{mi^2/2} z^i f_{[mi,r]}(rz) \, .
  &\equa{dualGFRemB}\cr
  }
$$
Setting $g(z)=\phi(r^{-m/2}z)$ and $g_i=[z^i]g(z)$, \dualGFRemB\ becomes
$$
  f(z)=1+f(z)\sum_{n\geq 1} g_n r^{mn^2/2} z^n f_{[mn,r]}(rz) \, .
$$
This expression is to be compared with (12.7) in Gessel (1980) which is (with
Gessel's $g_0$ set to $1$)
$$
  f(z)=1+\sum_{n\geq 1} g_n r^{mn^2/2} z^n f_{[mn,r]}(z) \, .
$$

In our setting, the analogue of Gessel's (1980) identity (12.7) is the second
assertion of Theorem \dualGF, which provides, 
with $H=\bigl(\lstar\,(1-\phi)\bigr)\,\widehat{\ }^{\,(p-1)}$,
$$
  f_{[k,r]}(z)=H(zr^k)/H(z) \, ,
$$
and the analogue of his (12.4) is that
$$
  \bigl(\lstar\,(1-\phi)\bigr)\,\widehat{\ }^{\,(p-1)}(z) 
  (\calT_{\phi,p})_{k,1/q}(z)
  = \bigl(\lstar\,(1-\phi)\bigr)\,\widehat{\ }^{\,(p-1)}(z/q^{k-1}) \, ,
$$
or, written differently, $H(z)f_{[k,r]}(z)=H(zr^k)$. This last identity
implies $[z^n]\bigl(H(z)f_{[k,r]}(z)\bigr)=r^{nk}[z^n]H(x)$, which is the 
analogue of Gessels's (1980) identity (12.11) in our setting.

\bigskip

\Remark Identity \dualGFEquation\ can be put in a fairly compact form using 
the operator $U_{p,\psi,q}$. Indeed, since $T_{0,1}=1$
(see Lemma \firstDualCoefficients),  the power series 
$\psi(z)=1-\calT_{\phi,p}(z)$ is of order $1$. Considering the right hand 
side of \dualGFEquation, we have
$$\eqalign{
  q^{-(p-1){i\choose 2}} z^i\prod_{0\leq j<(p-1)i} \calT_{\phi,p}(q^{-j}z)
  &{}= q^{-(p-1){i\choose 2}} \psi_{p,i,1/q}(z) \cr
  &{}= q^{(p-2)i/2}\bigl( U_{p,\psi,1/q}z^i\bigr) \, .
  \cr}
$$
Therefore, since $\calT_{\phi,p}=1-\psi$, \dualGFEquation\ can be rewritten as
$$\eqalign{
  -\psi(z)
  &{}=\sum_{i\geq 1} \phi_i q^{(p-2)i/2} U_{p,\psi,1/q}z^i \cr
  &{}=U_{p,\psi,1/q}\sum_{i\geq 1} \phi_i q^{(p-2)i/2} z^i \cr
  &{}=U_{p,\psi,1/q} \phi(q^{(p-2)/2}z) \, . \cr
  }
$$
Writing now $g(z)=-\phi(q^{p/2}z)$, we see that $\psi$ is the function $f$
that solves
$$
  f(z)=U_{p,f,1/q}g(z/q) \, .
$$
Up to substituting $q$ for $1/q$, this is exactly  Gessel's (1980) 
equation (12.1), which is $f(z)=U_{f,p,q}g(qz)$.

\bigskip

\Proof The proof has two parts corresponding to the two assertions.

\noindent{\it Part 1. Proof of the functional relation.}
In order to keep the subscripting within reason, we will 
write $t_n$ for $T_{n,n(p-1)+1}$.
As in Proposition \pAryDual, consider the Catalan power series
$P(z,t)=t-t\phi(t^{p-1}z)$, for which $t_n$ are the dual coefficients.
Those dual coefficients obey the recursion given by Theorem \Segner. To write
this recursion, we have, with $R(z,t)=\tilde P(z,t)t$,
$$
  R(z,t)
  = \phi(t^{p-1}z)/z  
  =\sum_{i\geq 1} \phi_i t^{(p-1)i} z^{i-1} \, .
$$
Thus, writing $R_{i,j}$ for the coefficient of $z^it^j$ in $R(z,t)$, we have
$$
  R_{i-1,(p-1)i}=\phi_i \, ,\quad i\geq 1\, ,
$$
and $R_{i,j}=0$ if $j\not=(p-1)(i+1)$ for nonnegative $i$. Writing $k_i$ 
for $\bigl(n_i,n_i(p-1)+1\bigr)$, recursion \SegnerRecursion\ takes the form
$$\displaylines{
  t_n
  = \One\{\, n=0\,\} +\sum_{i\geq 1} \phi_i \sum_{n_1,\ldots, n_{(p-1)i+1}\in\NN}
  t_{n_1}\ldots t_{n_{(p-1)i+1}}
  \hfill\cr\hfill
  q^{B_{(p-1)i+1}(k_1,\ldots, k_{(p-1)i+1})} \One\{\, n_1+\cdots+n_{(p-1)i+1}
  =n-i\,\}\, .
  ~\equa{dualGFA}\cr}
$$
Arguing as in the proof of Lemma \BnExpand, $B_m(k_1,\ldots,k_m)$ is
$$\displaylines{
  \sum_{1\leq i<j\leq m} \bigl( (p-1)n_i+1\bigr)n_j
  \hfill\cr\hfill
  = (p-1){\sum_{1\leq i\leq m} n_i\choose 2} 
    -(p-1)\sum_{1\leq i\leq m} {n_i\choose 2} 
    +\sum_{1\leq j\leq m} (j-1)n_j \, . \cr
  }
$$
Hence, \dualGFA\ yields
$$\displaylines{
  t_n = \One\{\, n=0\,\} 
  + \sum_{i\geq 1} \phi_i \sum_{n_1,\ldots,n_{(p-1)i+1}\in\NN}
  t_{n_1}\ldots t_{n_{(p-1)i+1}} 
  \hfill\cr\hfill
  q^{ (p-1){n-i\choose 2} -(p-1)\sum_{1\leq j\leq (p-1)i+1} {n_j\choose 2} 
      + \sum_{1\leq j\leq (p-1)i+1}(j-1)n_j }
  \hfill\cr\hfill
  \One\{\, n_1+\cdots+n_{(p-1)i+1}=n-i\,\} \, .\quad
  \equa{dualGFB}
}
$$
Set $\tilde t_n=q^{-(p-1){n\choose 2}} t_n$. We see that 
if $n_1+\cdots + n_{(p-1)i+1}=n-i$, then
$$\eqalign{
  {n-i\choose 2} -{n\choose 2} 
  &{}=-in + {i(i+1)\over 2}\cr
  &{}=-{i\choose 2} -i(n_1+\cdots+n_{(p-1)i+1}) \, .\cr}
$$
Multiplying both sides of \dualGFB\ by $z^n q^{-(p-1){n\choose 2}}$,
a quantity which is $1$ when $n=0$, yields
$$\eqalign{
  z^n \tilde t_n 
  &{}=\One\{\, n=0\,\} +\sum_{i\geq 1} \phi_i z^i q^{-(p-1){i\choose 2}}
    \hskip -4pt
    \sum_{n_1,\ldots, n_{(p-1)i+1}} \tilde t_{n_1}\ldots 
    \tilde t_{n_{(p-1)i+1}} \cr
  &\hskip 60pt {}\times q^{\sum_{1\leq j\leq (p-1)i+1}(j-1-(p-1)i) n_j} 
    z^{n_1+\cdots + n_{(p-1)i+1}} \cr\noalign{\vskip 3pt}
  &\hskip 60pt\qquad {}\times\One\{\, n_1+\cdots + n_{(p-1)i+1}=n-i\,\} \cr
  \noalign{\vskip 4pt}
  &{}=\One\{\, n=0\,\} +\sum_{i\geq 1} \phi_i z^i q^{-(p-1){i\choose 2}}
    \sum_{n_1,\ldots, n_{(p-1)i+1}} \Bigl( \prod_{1\leq j\leq (p-1)i+1}\cr
  &\hskip 30pt{}\tilde t_{n_j} q^{(j-1-(p-1)i)n_j} z^{n_j}\Bigr) 
    \One\{\, n_1+\cdots+n_{(p-1)i+1}=n-i\,\} \, . \cr}
$$
We sum over $n$ to obtain
$$
  \calT_{\phi,p}(z)
  = 1+\sum_{i\geq 1} \phi_i z^i q^{-(p-1){i\choose 2}} 
    \prod_{1\leq j\leq (p-1)i+1}\calT_{\phi,p}(q^{j-1-(p-1)i}z) \, ,
$$
that is, after substituting $j$ for $(p-1)i-j+1$ in the product, the 
first assertion of Theorem \dualGF.

\noindent{\it Part 2. Proof of the second assertion.} The proof is an
abstraction of some calculations in Prellberg and Brak (1995) ---~but see also
Gessel (1980, section 12) where simlar ideas are used. Considering the
functional relation, we see that it is nonlinear because of the product
term $\prod_{0\leq j\leq (p-1)i} \calT_{\phi,p}(q^{j-(p-1)i}z)$. As noticed
in Prellberg and Brak (1995) in some special cases, writing 
$\calT_{\phi,p}(z)$ as $g(z/q)/g(z)$ for some power series $g$ to be determined
simplifies this product considerably, and in fact linearizes the functional 
equation. Repeated substitutions in the relation $g(z)=g(z/q)/\calT_{\phi,q}(z)$
shows that
$g(z)=1/\prod_{j\geq 0} \calT_{\phi,q}(z/q^j)$ is this function; but this product
form is not useful beyond showing the existence and uniqueness of the power
series $g$. This change of power series leads to
$$\eqalign{
  \prod_{0\leq j\leq (p-1)i} \calT_{\phi,p}(q^{-j}z)
  &{}= \prod_{0\leq j\leq (p-1)i} {g(q^{-j-1}z)\over g(q^{-j}z)} \cr
  &{}= {g(q^{-(p-1)i-1}z)\over g(z)} \, .\cr
  }
$$
Thus, multiplying both sides of \dualGFEquation\ by $g(z)$, we obtain
$$
  g(z/q)=g(z)+\sum_{i\geq 1} \phi_i z^i q^{-(p-1){i\choose 2}} 
  g(q^{-(p-1)i-1}z) \, ,
  \eqno{\equa{dualGFC}}
$$
an equation which is linear in $g$. We write $g(z)=\sum_{n\in\NN}g_nz^n$ for the
power series expansion of $g$, and set $g_i=0$ if $i$ is negative. 
Considering the coefficient of $z^n$ in this equation, we obtain the identity
$$
  q^{-n} g_n = g_n +\sum_{i\geq 1} \phi_i 
    q^{-(p-1){i\choose 2}-(n-i)((p-1)i+1)} g_{n-i} \, .
   \eqno{\equa{dualGFD}}
$$
To understand where the Garsia roofing operator comes from in this 
calculation, note that given the complicated power of $q$ in right 
hand side of \dualGFD, it is quite natural
to make a change of variable, setting $g_n=q^{-\alpha_n}a_n$, and determine
the sequence $(\alpha_n)$ which allows us to rewrite \dualGFD\ in the
simplest manner. Making this change of variable, and multipliying
both sides of \dualGFD\ by $q^{\alpha_n}$ yields
$$
  q^{-n} a_n = a_n +\sum_{i\geq 1} \phi_i 
  q^{-(p-1){i\choose 2} -(n-i)((p-1)i+1)-\alpha_{n-i}+\alpha_n} a_{n-i} \, .
$$ 
The exponent of $q$ is
$$
  -{p-1\over 2}\bigl( n^2-(n-i)^2\bigr) +{p-1\over 2} i -n+i-\alpha_{n-i}
  +\alpha_n \, .
$$
Therefore, choosing $\alpha_n=(p-1)n(n-1)/2$ changes this exponent into $i-n$,
and yields the recursion
$$
  q^{-n} a_n=a_n+\sum_{i\geq 1} \phi_i q^{i-n} a_{n-i} \, ,
  \eqno{\equa{dualGFE}}
$$
with $a_i=0$ if $i$ is negative. Consider the generating function of $a_n$, 
that is, $A(z)=\sum_{n\geq 0} a_n z^n$.
Identity \dualGFE\ asserts that
$$\eqalign{
  [z^n] A(z/q)
  &{}= [z^n]A(z)+\sum_{i\geq 1} [z^i]\phi(z)[z^{n-i}]A(z/q) \cr
  &{}= [z^n] A(z)+[z^n]\bigl(\phi(z)A(z/q)\bigr)\, . \cr
  }
$$
Thus, $A(z/q)=A(z)+\phi(z)A(z/q)$, and therefore, 
$A(z)=\bigl(1-\phi(z)\bigr) A(z/q)$. By repeated 
substitutions, $A(z)=\prod_{i\geq 0} \bigl( 1-\phi(q^{-i}z)\bigr)$. This implies
$$
  g(z)=\sum_{n\geq 0} q^{-(p-1)n(n-1)/2} z^n [z^n]A(z) \, .
$$
Given the remark following Theorem \dualGF, we see that $A$ is $\Phi_{1/q}$
and that the function $g$ in this proof coincides with that given in the
remark. This proves the result.\hfill\qed

\bigskip

As a consequence of Theorem \dualGF, we can express the generating function
$\calT_{\phi,p}$ as a ratio of two constant terms. They involve
the Jacobi style theta function
$$
  \theta(u,q)=\sum_{n\in\NN} q^{n\choose 2} u^n \, .
$$

\Proposition{\label{GFConstantTerms}
  With the notation of Theorem \dualGF,
  $$
    \calT_{\phi,p}(z)
    = { [u^0] \biggl( \theta\Bigl({\ds z\over\ds uq},{\ds 1\over\ds q^{p-1}}
                           \Bigr)\Phi_{1/q}(u)\biggr)
        \over
        [u^0] \biggl( \theta\Bigl({\ds z\over\ds u},{\ds 1\over\ds q^{p-1}}
                           \Bigr)\Phi_{1/q}(u)\biggr)
      } \, . 
  $$
}

\Proof Note that 
$$
  [u^{-n}] \theta(z/u,1/q^{p-1}) = q^{-(p-1){n\choose 2}} z^n \, .
$$
Thus, the function $g$ in the remark following Theorem \dualGF\ can be
written as
$$\eqalign{
  g(z)
  &{}=\sum_{n\in\NN}\, [u^{-n}] \theta(z/u,1/q^{p-1})\,[u^n]\Phi_{1/q}(u) \cr
  &{}=[u^0] \bigl( \theta(z/u,1/q^{p-1})\Phi_{1/q}(u)\bigr) \, . \cr
  }
$$
The result follows from Theorem \dualGF.\hfill\qed

\bigskip


\section{\poorBold{$q$}-analogue of the Fuss-Catalan numbers}
A $p$-ary tree is a tree for which each parent has $p$ children. The 
Fuss-Catalan number $C_{p,n}$ counts the numbers of $p$-ary trees with $np+1$
nodes. From this definition, these numbers obey the recursion
$$
  C_{p,n}
  =\One\{\, n=0\,\} +\sum_{r_1,\ldots, r_p\in\NN} 
  \One\{\, r_1+\cdots+r_p=n-1\,\} C_{p,r_1}\cdots C_{p,r_p} \, .
  \eqno{\equa{FussCatalanDef}}
$$
Alternative definition in terms of lattice paths, disections of polygons,
or staircase tilings
exist, as shown for instance in Hilton and Perderson's (1991) and Heubach, Li
and Mansour (2008) papers. They also arise in some probabilistic problems
(Bajunaid, Cohen, Colonna, Singman, 2005; Liggett, 2000).
From the recursion \FussCatalanDef, the generating function
$$
  \widehat C_p(z)=\sum_{n\in\NN} C_{p,n} z^n
$$
solves the equation 
$$
  t^pz-t+1=0 \, .
  \eqno{\equa{FussCatalanEquation}}
$$ 
Lagrange inversion leads to the explicit formula
$$
  C_{p,n}={1\over (p-1)n+1} {pn\choose n} \, .
$$

We will now define some $q$-Fuss-Catalan numbers $(\calC_{p,n})_{n\in\NN}$ 
in the spirit of Carlitz's (1972) $q$-Catalan numbers. It is unclear from 
the definition that we will use that these numbers have a combinatorial 
interpretation in terms of lattice paths; while one may suspect, by analogy, 
some connection with $(p-1)$-Dyck paths, the results of 
the previous sections provide various interpretations for those numbers; and
we will see that one of those interpretations allows us to indeed derive a 
combinatorial one in terms of $(p-1)$-Dyck paths.

Given \FussCatalanEquation, consider the polynomial $P(z,t)=t-t^pz$. It is
a Catalan power series of the form $t-t\phi(t^{p-1}z)$ with $\phi(z)=z$. The
corresponding $p$-ary powers are
$$\eqalign{
  \phi_{p,n,q}(z)
  &{}=z^n\prod_{0\leq j<n(p-1)} \bigl( 1-\phi(q^jz)\bigr) \cr
  &{}=z^n (z,q)_{n(p-1)} \, . \cr
  }
$$
Proposition \pAryDual\ asserts that the corresponding dual coefficients
$(T_i)_{i\in\NN^2}$ can be identified with a sequence $(\calC_{p,n})_{n\in\NN}$
as follows. Define $(\calC_{p,n})_{n\in\NN}$ as the unique sequence such that
$$
  \sum_{n\in\NN} \calC_{p,n} z^n(z,q)_{n(p-1)+1} = 1 \, .
$$
We then have
$$
  T_i
  =\cases{ \calC_{p,n} & if $i=\bigl( n,n(p-1)+1\bigr)$ for some $n\in\NN$;\cr
           \noalign{\vskip 3pt}
           0 & otherwise.\cr}
$$

Alternatively, using Theorem \dualGF, the sequence $(\calC_{p,n})$ may be
defined through the generating function 
$$
  \calT(z)
  =\calT_{\phi,p}(z)
  =\sum_{n\in\NN}q^{-(p-1){n\choose 2}}\calC_{p,n}z^n
$$ 
by $\calT(0)=1$ and the functional relation
$$
  \calT(z)=1+z\prod_{0\leq j\leq p-1} \calT(q^{-j}z) \, .
  \eqno{\equa{FussCatalanFunctional}}
$$
Comparing this functional relation with display (12) in Bergeron (2012), we
see that $q^{-(p-1){n\choose 2}}\calC_{p,n}$ has a combinatorial interpretation
for it counts the area of $(p-1)$-Dyck paths. This section then yields new 
results concerning the area of these path, which are the analogue of known
results on the usual Dyck paths and the analogue of results on $q$-Catalan
numbers related to the usual Dyck paths.

Theorem \dualGF\ allows us to express this generating function with a new
form of $q$-Airy function. Indeed, we saw in the example following Theorem 
\dualGF\ that with $\phi(z)=z$ we have
$$
  [z^n]\Phi_{1/q}(z)=(-1)^n {q^{-{n\choose 2}}\over (1/q,1/q)_n} \, .
$$
We have, with the notation of the remark following Theorem \dualGF,
$$
  g(z)=\sum_{n\in\NN} {q^{-p{n\choose 2}} (-1)^n z^n\over (1/q,1/q)_n} \, . 
$$
This function $g$ is a new $q$-analogue of the Airy function, call it
$\Ai_{p,1/q}$. Then Theorem \dualGF\ asserts that
$$
  \calT(z)={\Ai_{p,1/q}(z/q)\over \Ai_{p,1/q}(z)} \, .
$$
Theorem \qEquation\ provides an alternative viewpoint. Let $(M,A)$ be two
$q$-commuting variables and consider the equation $A=T-T^pM$ with unknown
$T$. The unique power series in $(M,A)$ solving this equation is
$$
  T
  =\sum_{i\in\NN^2} T_i(M,A)^i
  =\sum_{n\in\NN} \calC_{p,n} M^n A^{n(p-1)+1}\, . 
$$

Finally, Theorem \Segner\ provides yet another definition. Indeed, with the 
notation of that theorem, $R(z,t)=t^{p-1}$ and therefore $R_{i,j}=\One\{\, i=0
\,;\, j=p-1\,\}$. Recursion \SegnerRecursion\ asserts that 
$$\displaylines{\quad
  \calC_{p,n}
  = \One\{\, n=0\,\} +\sum_{k_1,\ldots,k_p\in\NN} \calC_{p,k_1}\cdots
  \calC_{p,k_p}
  \hfill\cr\hfill
  q^{(p-1)\sum_{1\leq i<j\leq p} k_ik_j +\sum_{1\leq j\leq p} (j-1)k_j}
  \One\{\, k_1+\cdots+k_p=n-1\,\} \, ,
  \cr}
$$
which is a $q$-analogue of \FussCatalanDef\ and a $p$-ary analogue of
\qSegner.
Given this recursion, it appears that these $q$-Fuss-Catalan numbers are
different than those introduced by Haiman (1994), Garsia and Haiman (1996) 
or those of type A and B defined by Stump (2008; section 2.4).


\bigskip

\noindent{\bf References.}

{\leftskip=\parindent \parindent=-\parindent
 \par

G.\ Andrews (1975). Identities in combinatorics. II: A $q$-analog of the
Lagrange inversion theorem, {\sl Proc.\ Amer.\ Math.\ Soc.}, 53, 240--245.

G.\ Andrews, R.\ Askey, R.\ Roy (2000). {\sl Special Functions}, Cambridge.

I.\ Bahunaid, J.M.\ Cohen, F.\ Colonna, D.\ Singman (2005). Function series, 
Catalan numbers, and random walks on trees, {\sl The Amer.\ Math.\ Monthly},
112, 765--785.

F.\ Bergeron (2012). Combinatorics of $r$-Dyck paths, $r$-parking functions,
and the $r$-Tamari lattices, {\tt arxiv.1202.6269}.

L.\ Carlitz (1972). Sequences, paths, ballot numbers, {\sl Fibonacci Quart.},
10, 531--549.

Ph.\ Flajolet, R.\ Sedgewick (2009). {\sl Analytic Combinatorics}, Cambridge.

J.\ F\"urlinger, J.\ Hofbauer (1985). $q$-Catalan numbers, {\sl J.\ Comb.\
Th., A}, 248--264.

A.M.\ Garsia (1981). A $q$-analogue of the Lagrange inversion formula,
{\sl Houston J.\ Math.}, 7, 205--237.

A.M.\ Garsia, M.\ Haiman (1996). A remarkable $q,t$-Catalan sequence and
$q$-Lagrange inversion, {\sl J.\ Algebraic Comb.}, 5, 191-244.

I.\ Gessel (1980). A noncommutative generalization and $q$-analog of the 
Lagrange inversion formula, {\sl Trans.\ Amer.\ Math.\ Soc.}, 257, 455-482.

M.\ Haiman (1994). Conjectures on the quotient ring by diagonal invariants,
{\sl J.\ Algebraic\ Comb.}, 3, 17-76.

S.\ Heubach, N.Y.\ Li, T.\ Mansour (2008). Staircase tilings and 
$k$-Catalan structures, {\sl Discrete Mathematics}, 308, 5954-5953.

P.\ Hilton, J.\ Pedersen (1991). Catalan numbers, their generalizations, and
their uses, {\sl Math.\ Intelligencer}, 13, 64--75.

J.\ Hofbauer (1984).  A $q$-analog of the Lagrange expansion, {\sl Arch.\
Math.}, 42, 536--544.

M.E.H.\ Ismail (2009). {\sl Classical and Quantum Orthogonal Polynomials in One
Variable}, Cambridge.

T.\ Koshy (2009). {\it Catalan Numbers with Applications}, Oxford.

C.\ Krattenthaler (1984). A new $q$-Lagrange formula and some applications,
{\sl Proc.\ Amer.\ Math.\ Soc.}, 90, 338--344.

C.\ Krattenthaler (1988). Operator methods and Lagrange inversion: a unified
approach to Lagrange formulas, {\sl Trans.\ Amer.\ Math.\ Soc.}, 305, 431--465.

T.M.\ Liggett (2000). Monotonicity of conditional distributions and growth
models on trees, {\sl Ann.\ Probab.}, 28, 1645--1655.

J.\ McDonald (1995). Fiber polytopes and fractional power series, {\it J.\
Pure Appl.\ Alg.}, 104, 213--233.

T.\ Prellberg, R.\ Brak (1995). Critical Exponents from Non-Linear Functional 
Equations for Partially Directed Cluster Models, {\sl J.\ Stat.\ Phys.}, 
78, 701-730.

C.\ Stump (2008). { $q,t$-Fu\ss -Catalan numbers for finite reflection
groups}, PhD Thesis, Universit\"at Wien.

R.J.\ Walker (1950). {\it Algebraic Curves}, Princeton.

}

\medskip

\setbox1=\vbox{\halign{#\hfil &\hskip 40pt #\hfill\cr
  Ph.\ Barbe                  & W.P.\ McCormick\cr
  90 rue de Vaugirard         & Dept.\ of Statistics \cr
  75006 PARIS                 & University of Georgia \cr
  FRANCE                      & Athens, GA 30602  \ \ USA\cr
  philippe.barbe@math.cnrs.fr & bill@stat.uga.edu \cr}}
\dp1=0pt\ht1=0pt
\box1


\DoNotPrint{

\vfill\eject


\def\DATE{Sept.\ 15, 2012}

Consider the remark after Theorem \dualGF. Using Cauchy formula,
the function $g$ in that remark is
$$\eqalign{
  g(z)
  &{}=\sum_{n\in\NN} q^{-(p-1){n\choose 2}} z^n {1\over 2\pi i}\oint 
    {\Phi_{1/q}(\zeta)\over \zeta^{n+1}} \d\zeta \cr
  &{}={1\over 2\pi i} \oint \sum_{n\in\NN} q^{-(p-1){n\choose 2}} 
    \Bigl({z\over \zeta}\Bigr)^n {1\over \zeta} \Phi_{1/q}(\zeta) \d\zeta \, .
    \cr}
$$
Recall that Jacobi's theta function is (essentially, and for complex argument!)
$$
  \theta(r,y)=\sum_{n\in\NN} r^{-n^2} y^n \, .
$$
Since
$$
  (p-1){n\choose 2}={(p-1)\over 2} n^2- {p-1\over 2} n \, ,
$$
we obtain
$$
  g(z)
  ={1\over 2\pi i} \oint 
  \theta\Bigl(q^{(p-1)/2},{z\over\zeta}q^{(p-1)/2}\Bigr) {1\over \zeta}
  \Phi_{1/q}(\zeta) \d \zeta \, .
$$
So, to study the asymptotic behavior of $g$ as $q$ tends to $1$, hence, given
Theorem \dualGF, the asymtotic behavior of the generating function of
the dual coefficients of $q$-ary powers, hence, the asymptotic behavior
of the dual coefficients as $q$ tends to $1$, a possible first step is
to study the asymtotic behavior of the theta function and of $\Phi_{1/q}$.

That seems quite manageable. Write $q=e^\epsilon$, and we seek asymptotic
s $\epsilon$ tends to $0$. Some rough calculations, replacing
$$
  \log \Phi_{1/q}(\zeta)=\sum_{n\in\NN}\log\bigl( 1+\phi(e^{-n\epsilon}\zeta)\bigr)
$$
by the integral
$$
  \int_0^\infty \log\bigl(1+\phi(e^{-u\epsilon}\zeta)\bigr) \d u
$$
and making the change of variable $s=e^{-\epsilon u}\zeta$ yields
$$
  \log\Phi_{e^{-\epsilon}}(\zeta)
  \approx {1\over\epsilon} \int_0^\zeta {\log\bigl(1+\phi(s)\bigr)\over s}
  \d s \, .
$$
Now, for the theta function,
$$
  \theta(\epsilon,y)=\sum_{n\in\NN} e^{-n^2\epsilon+n\log y} \, ,
$$
the idea is to use the usual quadratic (Gaussian) approximation arround the 
term of the series which is maximal. 
The exponent is minimum when $n=(\log y)/2\epsilon$, and in this case is
about $(\log^2 y)/2\epsilon$. So we should have (considering only the terms
of exponential order in this rough calculation),
$$
  \theta(\epsilon,y)\approx  e^{(\log^2y)/2\epsilon} \, .
$$
This would give
$$
  g_{e^\epsilon}(z)
  \approx \oint_\zeta \exp\biggl( {1\over \epsilon} 
    \Bigl( {1\over 2} \log^2{z\over \zeta} 
           +\int_0^\zeta {\log\bigl(1+\phi(s)\bigr)\over s}\d s
    \Bigr)\biggr) \d \zeta\, .
$$
Now, presumably, we could us a saddle point, stationnary phase or 
something in that spirit, which leads to see for which $\zeta$ we can make 
the term in the exponential term as small as possible and deform the path
of integration accordingly. The function involved in the exponential is,
up to the factor $1/\epsilon$,
$$
  {1\over 2} \log^2 z - \log z\log \zeta +{1\over 2} \log^2 \zeta
  +\int_0^\zeta {\log\bigl( 1+\phi(s))\over s} \d s \, .
$$
Its partial derivative with respect to $\zeta$ is
$$
  {\log z\over \zeta} +{\log \zeta\over\zeta} 
  +{\log\bigl(1+\phi(\zeta)\bigr)\over\zeta}
  ={\log \Bigl( z\zeta \bigl( 1+\phi(\zeta)\bigr) \Bigr)\over \zeta} \, .
$$
This vanishes when 
$$
  z\zeta\bigl( 1+\phi(\zeta)\bigr)=1 \, .
  \eqno{\equa{saddle}}
$$ 
So that $\zeta=\zeta(x)$ should be the one poping out in the asymptotic.

It is interesting to see that this is in line what the case known of the
Catalan numbers: $\phi(x)=-x$. In that case the equation \saddle\ becomes
$z\zeta(1-\zeta)=1$, which gives
$$
  \zeta={1\pm\sqrt{1-4/z}\over 2} \, ,
$$
which is close enough to the Catalan generating function to suggests that
this may work.


\section{Garsia's powers}
In general, the dual coefficients are a a two dimensional array. However,
in example b) in sections 2 and 3, it was in fact one dimensional since
only the $T_{n,n+1}$, $n\in\NN$, could not vanish. More can be said on this
situation, and our next result suggests that it is related to $\tilde P(z,t)$
being a power series in $zt$, or $e_k$ being a Garsia power, whose definition
we recall.

\Definition{\label{GarsiaPower} 
  Let $\phi$ be a power series of order $1$. Its Garsia powers are the power
  series
  $$
    \phi_{k,q}(z)=q^{-{k\choose 2}} \prod_{0\leq j<k} \phi(q^jz) \, ,
  $$
  with $\phi_{0,q}=1$.
}

\bigskip

Note that the factor $q^{-{k\choose 2}}$ in this definition ensures that
$\phi_{k,q}$ is a power series of order $k$ with $[z^k]\phi_{k,q}(z)=1$.

The following is a characterization of Garsia's powers in the context of
$q$-Catalan bases and power series.

\Theorem{\label{GarsiaSetUp}
  Let $P$ be a Catalan power series with associated
  normalized $q$-Catalan basis $(e_k)$ and predual 
  basis $(\tilde e_i)_{i\in\NN^2}$. The following are equivalent:

  \medskip
  \noindent (i) there exists a power series $\phi$ of order $1$ such that
  $e_k=\phi_{k,q}$ for any $k\in\NN$;

  \medskip
  \noindent (ii) there exists some power series $(f_n)$ such that
  $$
    \tilde e_i(z,t)
    =z^{\langle i,\e_1-\e_2\rangle} f_{\langle i,\e_2\rangle}(zt)
  $$
  for any $i\in\NN^2$;

  \medskip
  \noindent (iii) $\tilde P(z,t)$ is a power series in $zt$.

  \medskip
  In this case,

  \medskip
  \noindent (A) $\phi(z)=zP(z,1)=P(1,z)$;

  \medskip
  \noindent (B) $f_n=e_n$;

  \medskip
  \noindent (C) $P(z,t)=\phi(zt)/z$;

  \medskip
  \noindent (D) the dual coefficients are determined by $T_i=0$ if
  $i\not\in\{\, (n,n+1)\,:\, n\in\NN\,\}$ and $\sum_{n\in\NN} T_{n,n+1}e_n(z)=z$.
}

\bigskip

Before proving this theorem, let us explain how it points to an interesting
conundrum. Consider a basis of Garsia's powers $e_k=\phi_{k,q}$. Assume that
we are interested in the numbers (power series in $q$) $\calC_n$ 
defined by
$$
  \sum_{n\geq 0} \calC_n e_n(z)=z \, .
  \eqno{\equa{conundrum}}
$$
In order to evaluate those $\calC_n$, one could determin the dual basis
of $(e_n)$, and, from identity \conundrum, obtain $\calC_n=[e_n]z$. Theorem
\dualBasis\ asserts that to calculate the dual basis it suffices to evaluate
the dual coefficients. But assertion (D) of Theorem \GarsiaSetUp\ asserts
that the dual coefficients are defined by the 
identity $\sum_{n\in\NN} T_{n,n+1}\e_n(z)=z$, which is precisely \conundrum. In 
other words, in this circumstances, Theorem \dualBasis\ is of no help 
whatsoever. A particular case occure for the $q$-Catalan numbers which are
obtained for $\phi(z)=z(1-z)$, and this may explain why finding a nice
explicit expression for these numbers has been such an elusive problem so
far.

\bigskip

\noindent{\bf Proof of Theorem \GarsiaSetUp.} {\it (i)${}\Rightarrow{}$(iii).} 
Since $\phi$ is of order $1$ and $e_k$ 
is 
normalized, there exists a power series $\psi$ such that 
$\phi(z)=z+z^2\phi(z)$. Since $(e_k)$ is a $q$-Catalan basis associated
to $P$,
$$
  {e_k(qz)\over e_k(z)} 
  = {\phi(q^kz)\over\phi(z)}
  = {P(z,q^k)\over P(z,1)} \, .
$$
Setting $t=q^k$, we obtain
$$
  P(z,t) = {P(z,1)\over\phi(z)} \phi(zt) \, .
$$
Simplifying both sides of this identity by $t$ and simplifying by $z$ the ratio
in the right hand side, we obtain
$$
  1+tz\tilde P(z,t)
  = {1+z\tilde P(z,1)\over 1+z\psi(z)} \bigl( 1+zt\psi(zt)\bigr) \, .
  \eqno{\equa{GarsiaSetUpA}}
$$
Introducing the power series $c(z)$ such that
$$
  {1+z\tilde P(z,1)\over 1+z\psi(z)}=1+zc(z) \, ,
$$
identity \GarsiaSetUpA\ becomes, after some simple manipulations,
$$
  t\tilde P(z,t)=c(z) + t\psi(z,t)+zt\psi(zt)c(z) \, .
$$
Applying $[t^0]$ on both sides yields $c(z)=0$, and rewriting the identity,
$$
  t\tilde P(z,t)=t\psi(zt) \, .
  \eqno{\equa{GarsiaSetUpB}}
$$
This proves (iii).

\noindent {\it (A) and (C).} Identity \GarsiaSetUpB\ implies
$$
  P(z,t)=t+zt^2 \psi(z,t) \, ,
$$
which, given how $\psi$ was defined gives both (A) and (C).

\noindent {\it (ii)${}\Rightarrow{}$(iii)} Assume that (ii) holds. Given
how $\tilde e_i$ is defined, this means
$$
  z^{\langle i,\e_1\rangle} \prod_{0\leq j<\langle i,\e_2\rangle}
  P(q^jz,t)
  = z^{\langle i,\e_1-\e_2\rangle} f_{\langle i,\e_2\rangle} (zt) \, .
$$
Then,
$$\eqalign{
  f_n(zt)
  &{}=z^n\prod_{0\leq j<n} P(q^jz,t) \cr
  &{}=(zt)^n \prod_{0\leq j<n} \bigl( 1+q^j zt\tilde P(q^jz,t)\bigr) \, . \cr
  }
$$
Therefore, for any nonnegative integer $n$,
$$
  {f_{n+1}(zt)\over ztf_n(zt)}
  = 1+q^n zt\tilde P(q^nz,t) \, .
$$
Using the variables $\tau=q^nzt$ and $y=q^nz$,
$$
  {q^n f_{n+1}(\tau/q^n)\over \tau f_n(\tau/q^n)}
  = 1+\tau\tilde P(y,\tau/y) \, .
$$
Since the left hand side of this identity does not depend on $y$, this means
that the function $y\mapsto \tilde P(y,\tau/y)$\note{are we using convergent power series? If not, should we restrict $\tilde P$ to be a polynomial to deal only with Laurent series when considering $\tilde P(y,\tau/y)$? Is there a better proof avoiding the division?} is constant. But
$$\eqalign{
  \tilde P(y,\tau/y)
  &{}=\sum_{i,j\in\NN}\tilde P_{i,j} y^i (\tau/y)^j \cr
  &{}=\sum_{n\in \ZZ} \sum_{i,j\in\NN} \tilde P_{i,j} \One\{\, i-j=n\,\}\tau^j
    \, . \cr
  }
$$
Thus, for any nonzero $n$, considering the term in $y^n$,
$$
  \sum_{j\in \NN} \tilde P_{j+n,j}\tau^j=0 \, .
$$
This forces $\tilde P_{j+n,j}=0$ for any nonnegative integer $j$ and any 
positive integer $n$. Then
$$
  \tilde P(z,t)=\sum_{j\geq 0} \tilde P_{j,j} z^jt^j
$$
and this proves (iii).

\noindent{\it (iii)${}\Rightarrow{}$(ii).} Under (iii), writing $\psi(z)$
for $\tilde P(z,1)$ or $\tilde P(1,z)$, we have $\tilde P(z,t)=\psi(zt)$
and
$$\eqalignno{
  \tilde e_i(z,t)
  &{}=z^{\langle i,\e_1\rangle} t^{\langle i,\e_2\rangle}
    \prod_{0\leq j<\langle i,\e_2\rangle} \bigl( 1+q^j zt\psi(q^jzt)\bigr) \cr
  &{}=z^{\langle i,\e_1-\e_2\rangle} (zt)^{\langle i,\e_2\rangle}
    \prod_{0\leq j<\langle i,\e_2\rangle} \bigl( 1+q^j zt\psi(q^jzt)\bigr) 
  &\equa{GarsiaSetUpC}\cr
  }
$$
and this proves (ii).

\noindent{\it (B).} Identity \GarsiaSetUpC\ also shows that
$$
  f_n(z)=z^n\prod_{0\leq j<n} \bigl( 1+q^jz\psi(q^jz)\bigr) \, .
$$
Setting $\phi(z)=P(1,z)=z+z^2\psi(z)$, this means 
$f_n=q^{-{n\choose 2}} \phi_{n,q}$. Thus, it remains to show that $\phi_{n,q}$
coincide with $e_n$, which we do now.

\noindent{\it (iii)${}\Rightarrow{}$(i).} Define $\psi(z)=\tilde P(z,1)$,
so that, since (iii) holds, $\tilde P(z,t)=P(zt,1)=\psi(zt)$. Then Lemma
\representation\ implies
$$
  e_k(z)
  =z^k \prod_{j\in\NN} {1+q^jz\psi(q^jz)\over 1+q^{j+k}z\psi(q^{j+k}z)} \, .
$$
Thus, setting $\phi(z)=P(1,z)$, we obtain $e_k=\phi_{k,q}$.

\noindent{\it (D).} Given (ii), the dual coefficients are defined by
$$
  \sum_{i\in\NN^2} T_i z^{\langle i,\e_1-\e_2\rangle} f_{\langle i,\e_2\rangle}
  (zt)=t \, .
$$
Multiplying both sides of this identity by $z$ and using the variable
$\tau=zt$ gives
$$
  \sum_{i\in\NN^2} T_i z^{\langle i,\e_1-\e_2\rangle+1} 
  f_{\langle i,\e_2\rangle} (\tau)=\tau \, .
$$
Considering this expression as a power series in $(z,\tau)$, we must
have $T_i=0$ if $\langle i,\e_1-\e_2\rangle+1$ does not vanish. If
$\langle ,\e_1-\e_2\rangle+1=0$, then 
$\langle i,\e_2\rangle=\langle i,\e_1\rangle+1$ and $i$ is of the form
$(n,n+1)$ for some nonnegative integer $n$. Thus, the dual coefficients
are defined by
$$
  \sum_{n\geq 0} T_{n,n+1}f_{n+1}(\tau)=\tau \, .
$$
Given assertion (B), this proves assertion (D).\hfill\qed


\section{Examples}
Theorems \dualBasis, \qEquation\ and \Segner\ can be applied in a variety
of ways, and the purpose of this section is to illustrate some of the
possibilities.

}


\DoNotPrint{

\def\DATE{Oct.\ 24, 2012}

\hfuzz=50pt
Reply to Identity concerning ${}^*(1+\phi)$, Oct. 21, 2012.

I could not follow how the third display is obtained. As you explain, 
let us see how the first term comes about. How can we obtain a term in $z^k$ 
in $(1+\phi)^*(z)$? We have
$$
  (1+\phi)^*(z)=\prod_{i\geq 0} \bigl(1+\phi(q^iz)\bigr) \, .
$$
A $z^k$ comes from some $z^j$ with $j\geq 1$ in some terms $(1+\phi^iz)$. Since
$j$ is at least $1$, we need to consider at most $k$ such terms. Writing
$m$ for the number of such terms we consider,
$$\displaylines{\qquad
  [z^k](1+\phi)^*(z)=\sum_{1\leq m\leq k} \sum_{0\leq i_1<\ldots<i_m}
  \sum_{j_1,\ldots ,j_m} \One\{\, j_1+\cdots+j_m=k\,\}
  \hfill\cr\hfill 
  [z^{j_1}]\bigl(1+\phi(q^{i_1}z)\bigr)\cdots 
  [z^{j_m}]\bigl(1+\phi(q^{i_m}z)\bigr) \, .
  \qquad\cr}
$$
So the coefficient we are looking for is
$$\displaylines{\quad
  \sum_{1\leq m\leq k} \sum_{0\leq i_1<\ldots<i_m}\sum_{j_1,\ldots ,j_m} 
  \One\{\, j_1+\cdots+j_m=k\,\}
  q^{i_1j_1}\phi_{j_1} \cdots q^{i_mj_m}\phi_{j_m}
  \hfill\cr\hfill
  {}=   \sum_{1\leq m\leq k} \sum_{0\leq i_1<\ldots<i_m}\sum_{j_1,\ldots ,j_m} 
  \One\{\, j_1+\cdots+j_m=k\,\}
  q^{i_1j_1+\cdots +i_mj_m}\phi_{j_1} \cdots \phi_{j_m}
  \, .\cr}
$$
We now permute the two inner sums. We then need to evaluate, 
for $j_1,\ldots,j_m$ fixed,
$$\eqalign{
  \sum_{0\leq i_1<\cdots <i_m} q^{j_1i_1+\cdots +j_mi_m} 
  &{}=\sum_{0\leq i_1\leq \ldots \leq i_m} 
    q^{j_1 i_1+j_2(i_2+1)+\cdots +j_m(i_m+m-1)}\cr
  &{}= q^{0j_1+1j_2+\cdots +(m-1)j_m} \sum_{0\leq i_1\leq\ldots\leq i_m} 
    q^{i_1j_1+i_2j_2+\cdots +i_mj_m} \, . \cr}
$$
It seems that 
$$
 \sum_{0\leq i_1\leq\ldots\leq i_m} 
    q^{i_1j_1+i_2j_2+\cdots +i_mj_m}
  ={1\over m!}\sum_{i_1,\ldots,i_m\geq 0} q^{i_1j_1+\cdots +i_mj_m}
$$
(it this correct?) which would give
$$
  {1\over (1-q^{j_1})\cdots (1-q^{j_m})} \, .
$$
So we would obtain
$$\eqalign{
  [z^k](1+\phi)^*(z)
  &{}=\sum_{1\leq m\leq k} \sum_{j_1,\ldots,j_m} 
    {\phi_{j_1}\over 1-q^{j_1}}\cdots {\phi_{j_m}\over 1-q^{j_m}} \cr
  &{}=\sum_{1\leq m\leq k} \Bigl( \sum_{j\geq 1} {\phi_j\over 1-q^j}\Bigr)^m\cr}
\, .
$$
And one could arrange a little the geometric sum over $m$.
I am not sure this calculation is correct, but it seems to me that your 
argument is off. But perhaps I am missing something.

}

\bye